\documentclass[english]{smfart}

\usepackage[all]{xy}
\usepackage[leqno]{amsmath}
\usepackage{amssymb}
\usepackage{smfthm}
\usepackage[english,francais]{babel}
\usepackage{graphics}
\usepackage{enumerate}

\theoremstyle{plain}
\newtheorem{theorem}[subsection]{Theorem}
\newtheorem{lemma}[subsection]{Lemma}
\newtheorem{propo}[subsection]{Proposition}

\newtheorem*{theoremnonnum}{Theorem}

\theoremstyle{definition}
\newtheorem{remark}[subsection]{Remark}
\newtheorem{definition}[subsection]{Definition}
\newcounter{eq}[subsection]
\newcommand\eq{\stepcounter{eq}}
\newcommand\num{\eq\tag{\thesubsection.\theeq}}

\newcommand{\X}{\widetilde X}
\newcommand{\Q}{{\widetilde{Q}}}
\newcommand{\N}{\mathbf{N}}
\newcommand{\C}{\mathbf{C}}
\newcommand{\R}{\mathbf{R}}
\newcommand{\M}{\mathcal M}
\renewcommand{\o}{\mathbf{O}}
\newcommand{\so}{\mathbf{SO}}
\renewcommand{\u}{\mathbf{U}}
\newcommand{\gl}{\mathbf{GL}}
\newcommand{\ssl}{\mathbf{SL}}
\newcommand{\lra}{\longrightarrow}
\newcommand\Sp{{\mathbf{Sp}}}
\newcommand{\sok}{{\boldsymbol{\mathfrak{so}}}}
\newcommand{\glk}{{\boldsymbol{\mathfrak{gl}}}}
\newcommand{\spk}{{\boldsymbol{\mathfrak{sp}}}}
\newcommand{\rs}{\mathrm{rs}}
\renewcommand\O{\mathcal{O}}
\renewcommand\P{\mathbf{P}}
\renewcommand\L{\mathcal{L}}

\newcommand{\Z}{\mathbf Z}
\newcommand{\G}{\mathbf G}

\DeclareMathOperator{\rk}{rk}
\DeclareMathOperator{\Ad}{Ad}
\DeclareMathOperator{\Aut}{Aut}
\DeclareMathOperator{\Ext}{Ext}
\DeclareMathOperator{\Hom}{Hom}

\DeclareMathOperator{\Quot}{Quot}
\DeclareMathOperator{\Spec}{Spec}
\DeclareMathOperator{\Mat}{Mat}
\DeclareMathOperator{\tr}{tr}

\newcommand\git{/\!\!/}

\title[Orthogonal and symplectic bundles on curves]{Orthogonal and symplectic bundles on curves and quiver representations}
\alttitle{Fibr\'es orthogonaux et symplectiques sur les courbes et
repr\'esentations de carquois}
\author{Olivier Serman}

\address{Laboratoire Paul Painlev\'e \\ UMR 8524 du CNRS \\ Universit\'e des Sciences et Technologies de Lille \\ Cit\'e scientifique \\ F-59 655 Villeneuve-d'Ascq Cedex \\ France}

\email{Olivier.Serman@math.univ-lille1.fr}

\begin{document}

\frontmatter

\selectlanguage{english}
\begin{abstract}
 We show how quiver representations and their invariant theory naturally arise in the study of some moduli spaces parametrizing bundles defined on an algebraic curve, and how they lead to fine results regarding the geometry of these spaces.
\end{abstract}

\begin{altabstract}
On montre comment la th\'eorie des repr\'esentations de carquois
appara\^it naturellement lors de l'\'etude des espaces de modules
de fibr\'es principaux d\'efinis sur une courbe alg\'ebrique, et
comment elle permet d'analyser la g\'eom\'etrie de ces
vari\'et\'es.
\end{altabstract}



\altkeywords{Espaces de modules, fibr\'es orthogonaux et
symplectiques sur les courbes, repr\'esentations de carquois}

\maketitle

\mainmatter

\section*{Introduction}

Let $X$ be a smooth projective curve defined over an
algebraically closed field $k$ of characteristic $0$. It is a
natural question to try to find an algebraic variety which
parametrizes objects of some given kind defined on the curve $X$.

A first example is provided by the study of line bundles of degree
$0$ on $X$. It has been known essentially since Abel and Jacobi
that there is actually an abelian variety, the Jacobian variety
$J_X$, which parametrizes line bundles of degree $0$ on $X$. We
know a great deal about this variety, whose geometry is closely
related to the geometry of $X$.

Weil's suggestion in \cite{weil} that vector bundles (which appear
in his paper as ``$\gl_r$-divisors'') should provide a relevant
non-abelian analogue of this situation opened the way to a large
field of investigations, which led to the construction in the
1960's of the moduli spaces of \textit{semi-stable} vector bundles
of given rank and degree on $X$, achieved mainly by Mumford,
Narasimhan and Seshadri. Ramanathan then extended this
construction to prove the existence of moduli spaces for
semi-stable principal $G$-bundles on $X$ for any connected
reductive group $G$.

These varieties, which will be denoted by $\M_G$ in this paper,
have been intensively investigated since their construction,
especially for $G=\gl_r$. They have more recently drawn new
attention for the fundamental role they appeared to play in
various subjects, like Conformal Field Theory or Langland's
geometric correspondence.

In these notes we consider the following question:

\vspace{0.2cm}

\noindent If $H \lra G$ is a morphism between two reductive
groups, what can we say about the induced morphism $\M_H \lra
\M_G$ between moduli spaces?

\vspace{0.2cm}

This is a frequently encountered situation. For example, choosing
for $H$ a maximal torus $T \simeq (\G_m)^l$ contained in $G$ gives
a morphism from the moduli space $\M^0_T$ of topologically trivial
$T$-bundles (which is isomorphic to $(J_X)^l$) to the variety
$\M_G$. When $X$ is the projective line $\P^1$, we know from
\cite{grothendieck} that any principal $G$-bundle on $\P^1$ comes
from a principal $T$-bundle. If $X$ is an elliptic curve,
\cite{laszloell} shows that the morphism $\M^0_T \lra \M_G$ is a
finite morphism from $\M^0_T \simeq X^l$ onto the connected
component of $\M_G$ consisting of topologically trivial
semi-stable $G$-bundles. For higher genus curves, let us just say
that the morphism $\M^0_{\G_m}=J_X \lra \M_{\ssl_2}$, which sends
a line bundle $L$ to the vector bundle $L \oplus L^{-1}$, gives a
beautiful way to investigate the geometry of the moduli spaces of
semi-stable rank $2$ vector bundles on $X$ (see \cite{beauville2}).

We study here the case of the classical groups $H=\o_r$ and
$\Sp_{2r}$, naturally embedded in the general linear group. The
moduli variety $\M_{\o_r}$ then parametrizes semi-stable
orthogonal bundles $(E,q)$ of rank $r$ on $X$, and the morphism
$\M_{\o_r} \lra \M_{\gl_r}$ just forgets the quadratic form $q$.
In the same way, $\M_{\Sp_{2r}}$ parametrizes semi-stable
symplectic bundles, and $\M_{\Sp_{2r}}\lra \M_{\ssl_{2r}}$ forgets
the symplectic form. We will also consider $\so_r$-bundles, which
are {oriented} orthogonal bundles $(E,q,\omega)$, that is
orthogonal bundles $(E,q)$ together with an orientation, which is
defined as a section $\omega$ in $H^0(X,\O_X)$ satisfying $\tilde
q (\omega)=1$ (where $\tilde q$ is the quadratic form on $\det E
\simeq \O_X$ induced by $q$).

We have shown in \cite{comp} that the forgetful morphisms
$$\M_{\o_r} \lra \M_{\gl_r} \text{ \ \ and \ \ } \M_{\Sp_{2r}} \lra \M_{\ssl_{2r}}$$
\noindent are both closed immersions. In other words, these
morphisms identify the varieties of semi-stable orthogonal and
symplectic bundles with closed subschemes of the variety of all
vector bundles. Note that this means that the images in
$\M_{\gl_r}$ of these two forgetful morphisms are normal
subschemes. The proof involves an infinitesimal study of these
varieties, which naturally leads to some considerations coming
from representation theory of quivers (for example, we use the
fact that $\M_{\gl_r}$ is locally isomorphic to the variety
parametrizing semi-simple representations of a given quiver). We
present in Section \ref{imm} a proof of this result which
simplifies a little the one given in \cite{comp}.

The moduli spaces $\M_G$ are in general not regular (nor even
locally factorial), and a basic question is to describe their
singular locus and the nature of the singularities. If $X$ has
genus $g \geqslant 2$, the singular locus of $\M_{\ssl_r}$ has a
nice description, which has been known for long (see \cite{NR}): a
semi-stable vector bundle defines a smooth point in $\M_{\ssl_r}$
if and only if it is a \textit{stable} vector bundle, except when
$r=2$ and $g=2$ (in this very particular case, $\M_{\ssl_2}$ is
isomorphic to $\P^3$). For $G$-bundles one has to consider
\textit{regularly stable} bundles, which are stable $G$-bundle $P$
whose automorphism group $\Aut_G(P)$ is equal to the center $Z(G)$
of $G$. Such a bundle defines a smooth point in $\M_G$, and one
can expect the converse to hold, barring some particular cases.

We solve this question for classical groups. Using Schwarz's
classification \cite{schwarz} of {coregular representations}, we
prove in Section \ref{sing} that the smooth locus of $\M_{\so_r}$
is exactly the regularly stable locus, except when $X$ has genus
$2$ and $r=3$ or $4$. For symplectic bundles we prove that the
smooth locus of $\M_{\Sp_{2r}}$ is exactly the set of regularly
stable symplectic bundles (for $r \geqslant 2$). This proof, which
requires a precise description of bundles associated to points of
the moduli spaces, cannot be extended to another group $G$ without
a good understanding of the nature of these bundles.

\

\noindent\textsc{Acknowledgements.} It is a pleasure to thank
Michel Brion for being responsible of such an enjoyable and
successful Summer school, and for having let me take part in this
event. In addition to the occasion of spending two amazing weeks
in Grenoble, it was a unique opportunity to add to \cite{comp} a
new part which could not have found a better place to appear.

\section{The moduli spaces $\M_G$}

\subsection{} Let $X$ be a smooth projective irreducible curve of genus $g \geqslant 1$,
defined over an algebraically closed field of characteristic $0$.

We can associate to $X$ its \textit{Jacobian variety} $J_X$, which
parametrizes line bundles of degree $0$ on the curve. It is a
projective variety, whose closed points correspond bijectively to
isomorphism classes of degree $0$ line bundles on $X$. Moreover,
$J_X$ has the following \textit{moduli property}:
\begin{itemize}
 \item if $\mathcal L$ is a family of degree $0$ line bundles on $X$ parametrized by a scheme $T$, the classifying map $\varphi$ which maps a point $t \in T$ on the point in $J_X$ associated to the line bundle $\mathcal L_t$ defines a morphism $\varphi \colon T \lra J_X$,
 \item $J_X$ is ``universal'' for this property.
\end{itemize}
We should also mention here that $J_X$ comes with a (non-unique)
\textit{Poincar\'e bundle} $\mathcal P$ on the product $J_X \times
X$. It is a line bundle on $J_X \times X$, whose restriction
$\mathcal P_a$ to $\{a\} \times X$ is exactly the line bundle
associated to the point $a \in J_X$.

The Jacobian variety inherits many geometric properties from its moduli interpretation: let us just note here that it is an abelian variety which naturally carries a principal polarization. This extra data allows to describe sections of line bundles on $J_X$ in terms of \textit{theta functions}. This analytical interpretation of geometric objects defined on $J_X$ provides a powerful tool to investigate the beautiful relations between the curve and its Jacobian.

\subsection{} It has thus been natural to look for some possible generalizations of this situation. To do this, we can remark that line bundles are exactly principal $\G_m$-bundles. Replacing the multiplicative group $\G_m$ by any reductive group $G$ leads to the consideration of \textit{principal $G$-bundles} on $X$.

When $G$ is the linear group $\gl_r$, they are vector bundles on
$X$. Topologically, vector bundles on the curve $X$ are classified
by their rank $r$ and degree $d$, and the natural question is to
find an algebraic variety whose points correspond to isomorphism
classes of vector bundles on $X$ of fixed rank and degree. The
idea that such varieties parametrizing vector bundles should exist
and give the desired non-abelian generalization of the Jacobian
variety goes back to Weil (see \cite{weil}). However, the
situation cannot be as simple as it is for line bundles. Indeed,
the collection $\mathcal V_{r,d}$ of all vector bundles of rank
$r$ and degree $d$ on $X$ is not \textit{bounded}: we cannot find
any family of vector bundles parametrized by a scheme $T$ such
that every vector bundle in $\mathcal V_{r,d}$ appears in this
family. So we need to exclude some bundles in order to have a
chance to get a variety enjoying a relevant moduli property.

As we have said in the introduction, the construction of these
moduli spaces of vector bundles on $X$ has been carried out in the
1960's, mainly by Mumford and by Narasimhan and Seshadri. They
happened to show that one has to restrict to \textit{semi-stable}
bundles to obtain a reasonable moduli variety. This notion was
introduced first by Mumford in \cite{mumfordICM} in the light of
Geometric Invariant Theory.

Let us define the \textit{slope} of a vector bundle $E$ as the
ratio $\mu(E)=\deg(E)/\rk(E)$.

\begin{definition} \label{semistable}
A vector bundle $E$ on $X$ is said to be \textit{stable} (resp. \textit{semi-stable}) if we have, for any proper subbundle $F \subset E$, the slope inequality
$$\mu(F)<\mu(E) \text{\ \ \  (resp. } \mu(F) \leqslant \mu(E) \text{).}$$
\end{definition}

We will mainly be concerned in the following with degree $0$ vector bundles. In this case, saying that a bundle is stable just means that it does not contain any subbundle of degree $\geqslant 0$.

Mumford's GIT allowed him to provide the set of isomorphism classes of stable bundles of given rank and degree with the structure of a quasi-projective variety.

\begin{theorem}[Mumford]\label{mum}
 There exists a coarse moduli scheme $\mathcal U^\mathrm{st}_X(r,d)$ for stable vector bundles of rank $r$ and degree $d$ on $X$. Its points correspond bijectively to isomorphism classes of stable bundles of rank $r$ and degree $d$.
\end{theorem}

\noindent This result precisely means that, if $F_{X,r,d}^\mathrm{st}$ denotes the moduli functor which associates to a scheme $T$ the set of isomorphism classes of families of stable vector bundles of rank $r$ and degree $d$ on $X$ parametrized by $T$,
\begin{enumerate}[(i)]
\item there is a natural transformation $\varphi : F^\mathrm{st}_{X,r,d} \lra \Hom(-,{\mathcal U^\mathrm{st}_X(r,d)})$ such that any natural transformation $F^\mathrm{st}_{X,r,d} \lra \Hom(-,N)$ factors through a unique morphism $\mathcal U^\mathrm{st}_X(r,d) \lra N$,
\item the set of closed points of $\mathcal U^\mathrm{st}_X(r,d)$ is identified (via $\varphi$) to the set $F^\mathrm{st}_{X,r,d}(\Spec k)$ of isomorphism classes of stable vector bundles of rank $r$ and degree $d$.
\end{enumerate}

\noindent (Of course, the natural transformation $\varphi$ associates to a family $\mathcal F$ of stable bundles parametrized by $T$ the corresponding classifying morphism $t \in T \longmapsto  \mathcal F_t \in \mathcal U_X^\mathrm{st}(r,d)$.)

In particular, once we agree to exclude non stable bundles, we
obtain a collection of vector bundles which carries a natural
algebraic structure. Hopefully, those bundles that we have to
forget form a very small class inside the set of all vector
bundles, at least when $X$ has genus $g \geqslant 2$. Indeed,
stability (as well as semi-stability) is an open condition: if
$\mathcal F$ is a family of vector bundles of rank $r$ on $X$
parametrized by $T$, the stable locus $T^\mathrm{st}=\{t \in T |
\mathcal F_t \text{ is stable}\}$ is open in $T$ (see also Remark
\ref{codim}).

\subsection{} Almost simultaneously, Narasimhan and Seshadri found the same notion of stability, but from a completely different approach inspired by Weil's seminal paper \cite{weil}. The key observation is that the Jacobian $J_X$ of a complex curve $X$ is a complex torus, which can be topologically identified with the space $\Hom (\pi_1(X), S^1)$ of all $1$-dimensional unitary representations of the fundamental group $\pi_1(X)$ of $X$. This transcendental correspondence between unitary characters of $\pi_1(X)$ and line bundles on $X$ is obtained as follows: if $\widetilde X$ is a universal covering of $X$, we associate to a character $\rho \colon \pi_1(X) \to S^1$ the line bundle $L_\rho$ on $X$ defined as the quotient $\widetilde X \times^\rho \C$ of the trivial bundle $\widetilde X \times \C$ by the action of the fundamental group given by $(x,\lambda)\cdot \gamma = (x \cdot \gamma, \rho(\gamma)^{-1} \lambda)$ for all $\gamma \in \pi_1(X)$ (in other words, $L_\rho$ is the 
 $\pi_1(X)$-invariant subbundle of the direct image of the trivial line bundle on $\widetilde X$). Moreover, this bijection becomes an actual isomorphism for the complex structure induced on $\Hom(\pi_1(X),S^1) \simeq H^1(X,\R)/H^1(X,\Z)$ by the natural isomorphism between $H^1(X,\R)$ and $H^1(X,\O_X)$ given by Hodge theory.

As Weil suggested, \textit{unitary representations of the
fundamental group} of $X$ had to play a prominent role in the
study of vector bundles, if only because two unitary
$r$-dimensional representations $\rho_1$ and $\rho_2$ of
$\pi_1(X)$ give isomorphic vector bundles $E_{\rho_i}=\widetilde X
\times^{\rho_i} \C^r$ if and only if they are equivalent (this
does not hold any longer for arbitrary linear representations,
which ultimately led to the notion of \textit{Higgs bundles}).
However, there is a major difference for higher rank vector
bundles: this construction does not allow to obtain every degree
$0$ vector bundle. The main result in \cite{weil} states that a
vector bundle on $X$ can be defined by a linear representation
$\rho \colon \pi_1(X) \lra \gl_r$ if and only if it is a direct
sum of indecomposable degree $0$ vector bundles, so that we
already miss the (non semi-stable) rank $2$ vector bundles $L
\oplus L^{-1}$ with $\deg L \geqslant 1$. And if we consider only
unitary representations, we have to exclude more bundles.

The fundamental result \cite{NarasimhanSeshadri}, which is
``already implicit in the classical paper of A. Weil'', states
that, if $X$ has genus at least $2$, stable vector bundles of rank
$r$ and degree $0$ correspond exactly to equivalence classes of
irreducible unitary representations $\pi_1(X) \lra \u_r$ of the
fundamental group. As a consequence, Theorem \ref{mum} also shows
that the set of equivalence classes of irreducible unitary
representations of $\pi_1(X)$ has a natural complex structure,
depending on that of $X$.

\subsection{} Seshadri went further and constructed a compactification of this variety by considering \textit{unitary bundles}, i.e. bundles associated to any unitary representation. These bundles, which are also called \textit{polystable bundles}, are exactly direct sums of stable bundles of degree $0$ (more generally, we say that a semi-stable bundle of arbitrary degree is polystable if it splits as the direct sum of stable bundles). Using Mumford's theory, he obtained a projective variety $\mathcal U_X(r,0)$ which parametrizes isomorphism classes of polystable bundles of rank $r$ and degree $0$ on $X$, and contains $\mathcal U_X^\mathrm{st}(r,0)$ as an open subscheme.

\label{JH} However, no moduli property can be formulated in terms
of polystable bundles. Indeed polystability behaves very badly in
family; it is not even an open condition. There is in fact a more
natural way to think about the variety $\mathcal U_X(r,0)$, based
on the following relation between polystable and semi-stable
vector bundles. The crucial fact is that Jordan-H\"older theorem
holds in the category of all semi-stable vector bundles of degree
$0$ on $X$, so that we can associate to any such bundle $E$ the
Jordan-H\"older graded object $\mathrm{gr}\, E$ (sometimes called
\textit{semisimplification} of $E$), which is defined as the
direct sum $\mathrm{gr}\, E=\bigoplus F_i/F_{i-1}$ of the stable
subquotients given by any Jordan-H\"older composition series
$0=F_0 \subset F_1 \subset \cdots \subset F_l=E$ for $E$. We say
that two semi-stable vector bundles are $S$-equivalent if the
associated graded objects are isomorphic. The point is that
$S$-equivalence classes of degree $0$ vector bundles coincide with
isomorphism classes of polystable bundles: the equivalence class
of a vector bundle $E$ is characterized by the isomorphism class
of the corresponding graded object, which is a polystable bundle.

Seshadri proved that the classifying map $t \in T \longmapsto \mathrm{gr}\, \mathcal E_t$ associated to any family $\mathcal E$ of semi-stable vector bundles of rank $r$ and degree $0$ on $X$ parametrized by a variety $T$ defines a morphism $T \lra \mathcal U_X(r,0)$, and that the variety $\mathcal U_X(r,0)$ is in fact a coarse moduli space for semi-stable vector bundles of rank $r$ and degree $0$, whose closed points correspond to $S$-equivalence classes of vector bundles.

In arbitrary degree the corresponding result also holds:

\begin{theorem}[Seshadri]
There exists a projective variety $\mathcal U_X(r,d)$ which is a coarse moduli scheme for semi-stable bundles of rank $r$ and degree $d$ on $X$. Its closed points correspond bijectively to $S$-equivalence classes of semi-stable vector bundles, or, equivalently, to isomorphism classes of \textit{polystable} vector bundles of rank $r$ and degree $d$ on $X$.
\end{theorem}

It is a normal irreducible projective variety. Moreover, when $X$ has genus $g \geqslant 2$, it has dimension
$r^2(g-1)+1$, and contains as
a dense open subscheme the moduli variety $\mathcal
U^\mathrm{st}_X(r,d)$ of stable vector bundles. Although it is
commonly denoted by $\mathcal{U}_X(r,d)$, we will preferably use
here the notation $\M_{\gl_r}^d$, which keeps track of the
identification between rank $r$ vector bundles and principal
$\gl_r$-bundles.

\begin{remark}\label{codim}
 If $X$ is a curve of genus $g \geqslant 2$, it is not difficult to show that, unless $g=2$ and $r=2$, the strictly semi-stable locus $\mathcal U_X(r,d) \setminus \mathcal U_X^\mathrm{st}(r,d)$ is a closed subscheme of codimension at least $2$, which means that stable bundles represent a very large part of the set of all semi-stable bundles. In the same way, semi-stable bundles form a very large class inside the collection of all vector bundles. More precisely, we can show using Harder-Narasimhan filtrations that, if $\mathcal F$ is a family of vector bundles on $X$ parametrized by a smooth scheme $T$ such that the Kodaira-Spencer infinitesimal deformation map $T_t \lra \Ext^1(\mathcal F_t,\mathcal F_t)$ is everywhere surjective, then the complement $T \setminus T^\mathrm{ss}$ of the semi-stable locus $T^\mathrm{ss}$ has codimension at least $2$ (see \cite[4.IV]{LPV}).
\end{remark}

\subsection{} Building up on these ideas, Ramanathan considered in his thesis \cite{ramanathanthesis} the case of principal $G$-bundles on a curve $X$ for any complex connected reductive group $G$. Topologically, principal bundles with connected structure group $G$ on $X$ are classified by their topological type which is a discrete invariant belonging to $H^2(X,\pi_1(G)) \simeq \pi_1(G)$. Ramanathan's aim was to construct coarse moduli schemes for $G$-bundles on $X$ of a given topological type $\delta \in \pi_1(G)$. We have of course to restrict ourselves to a certain class of $G$-bundles.

The first step is to define semi-stability for principal
$G$-bundles. It is done by considering reductions of structure
group to parabolic subgroups of $G$. Here we need to recall a few
definitions involving principal bundles (see \cite{serre}). If $P$
is a $G$-bundle on $X$ and $F$ a quasi-projective variety acted
upon by $G$, the associated fiber bundle $P(F)$ (also denoted by
$P \times^G F$) is the fiber bundle defined as the quotient $(P
\times F) / G$, where $G$ acts on $P \times F$ by $(p,f) \cdot g =
(p \cdot g, g^{-1} \cdot f)$. If $\rho \colon G \to G'$ is a group
morphism, $P(G')$ is a principal $G'$-bundle: it is called
\textit{extension of structure group} of $P$ from $G$ to $G'$ and
sometimes denoted by $\rho_\ast P$. Conversely, if $P'$ is a
$G'$-bundle, we call \textit{reduction of structure group} of $P'$
to $G$ a pair $(P,\alpha)$ consisting of a $G$-bundle $P$ on $X$
together with an isomorphism $\alpha \colon P(G')
\buildrel\sim\over\lra P'$ between the associated bundle
$\rho_\ast P=P(G')$ and $P'$.

Reductions of structure group of a $G'$-bundle $P'$ correspond to
sections $\sigma$ of the fiber bundle $P' / G = P' \times^{G'}G'/G
\to X$ as follows: if $\sigma \colon X \to P'/G$ is such a
section, the pull-back via $\sigma$ of the $G$-bundle $P' \to
P'/G$ defines a $G$-bundle $\sigma^\ast P'$ on $X$ whose extension
$\sigma^\ast(P') \times^G G'$ is naturally isomorphic to $P'$.
Moreover, two sections give isomorphic $G$-bundles if and only if
they differ by an automorphism of $P'$.

\begin{definition}[Ramanathan] \label{Gss}
A $G$-bundle $P$ on $X$ is \textit{stable} (resp.
\textit{semi-stable}) if, for every parabolic subgroup $\Pi \subset G$, for every non trivial dominant character $\chi$ of $\Pi$, and for every $\Pi$-bundle $Q$ defining a reduction of structure group of $P$ to $\Pi$, the line bundle $\chi_\ast Q$ has degree $\deg (\chi_\ast Q) <  0$ (resp. $\leqslant 0$).
\end{definition}

This seemingly technical definition gives back for $G=\gl_r$ the classical definition \ref{semistable}. Moreover, in characteristic $0$, a $G$-bundle $P$ is semi-stable if and only if its adjoint vector bundle $\Ad(P)=P \times^G \mathfrak g$ is. (In positive characteristic, we need to introduce strongly semi-stable bundles to get an analogous result.)

\subsection{} \label{unitary}Then we need to know how $S$-equivalence has to be generalized. We have to recall the following facts (and refer to \cite{ramanathanthesis} for details). Each equivalence class defines a Levi subgroup $L \subset G$ and a stable $L$-bundle $Q$ such that the associated bundle $Q(G)$ belongs to the given class. Moreover, the $G$-bundle $Q(G)$ is uniquely defined, up to isomorphism, by its equivalence class. This bundle is the analogue of the Jordan-H\"older graded object characterizing $S$-equivalence classes for vector bundles. Such bundles are called \textit{unitary $G$-bundles}.

It should be noted that the theorem of Narasimhan and Seshadri
remains true in this context (whence the terminology of unitary
bundles): if $K \subset G$ denotes a maximal compact subgroup of a
connected semisimple group $G$, then any morphism $\pi_1(X) \lra
K$ defines a unitary $G$-bundle on $X$, and we get in this way a
bijection between conjugacy classes of representations $\pi_1(X)
\lra K$ and isomorphism classes of unitary $G$-bundles (see
\cite{ramanathan75} for the corresponding statement for connected
reductive groups).

Now we can recall the main result of \cite{ramanathanthesis}.

\begin{theorem}[Ramanathan]
 Let $G$ be a complex connected reductive group and $\delta \in \pi_1(G)$. There exists a coarse moduli scheme $\M_G^\delta$ for semi-stable principal $G$-bundles on $X$ of topological type $\delta$. It is an irreducible normal projective variety, whose points correspond bijectively to $S$-equivalence classes of semi-stable $G$-bundles.
\end{theorem}

\subsection{} \label{construction}
Let us briefly recall the main lines of the construction of $\M_G$
for a semisimple group $G$ (following \cite{BLS}).

We fix a faithful representation $\rho \colon G \lra \ssl_r$, an ample line bundle $\O_X(1)$ on $X$, and an integer $M$ such that, for every semi-stable $G$-bundle $P$, the rank $r$ vector bundle
$P(\ssl_r) \otimes \O_X(M)$ is generated by its global sections and satisfies
$H^1(X,  P(\ssl_r) \otimes \O_X(M))=0$. Let us
consider the functor $\underline{R}_G$ which associates to a scheme
$S$ the set of isomorphism classes of pairs $(\mathcal P,\alpha)$ consisting
of a $G$-bundle $\mathcal P$ over $S\times X$ with semi-stable
fibers together with an isomorphism $\alpha \colon \O_S^\chi
\buildrel \sim \over \lra {p_S}_\ast\left( \mathcal P(\ssl_r) \otimes
p_X^\ast\O_X(M)\right)$ (where $\chi=r(M+1-g)$, and $p_X$ and $p_S$ denote the projections from $S \times X$ onto $X$ and $S$). This functor, which is
introduced to relate $G$-bundles to vector bundles, is
representable by a \textit{smooth} scheme $R_G$, which will be
referred to as a \textit{parameter scheme}. The functor $\underline R_{\ssl_r}$
is indeed representable by a locally closed subscheme of the Hilbert
scheme ${\Quot}_{\O_X^\chi}^{r,rM}$. If $(\mathcal
U,u)$ denotes the universal pair on $R_{\ssl_r} \times X$, we can see that $\underline R_G$ is exactly the functor of global sections of $\mathcal U/G$. This functor is representable by a
smooth scheme $R_G$, which is affine over $R_{\ssl_r}$.

We know from Simpson's construction that the moduli scheme $\M_{\ssl_r}$ is the (good) quotient
$R_{\ssl_r}\git \Gamma$  of the parameter scheme $R_{\ssl_r}$ by the natural
action of $\Gamma=\gl_\chi$ (for sufficiently high $M$). The point is that $R_{\ssl_r}$ is exactly the open subset of semi-stable points for the action of $\Gamma$ on a closed subscheme of ${\Quot}_{\O_X^\chi}^{r,rM}$.

The parameter scheme $R_G$ also carries a natural action of
$\Gamma$, for which the structural morphism $R_G \lra R_{\ssl_r}$
is $\Gamma$-equivariant. A good quotient $R_G \git \Gamma$, if it
exists, provides the desired coarse moduli space for semi-stable
$G$-bundles. According to \cite[Lemma 5.1]{ramanathanthesis}, its
existence follows from the one of $R_{\ssl_r}\git \Gamma$.

This construction can be adapted to more general cases, and in
particular to $\gl_r$- and $\o_r$-bundles: we find that
$\M_{\gl_r}^0$ is the good quotient $R^0_{\gl_r} \git \Gamma$ of a
smooth parameter scheme $R^0_{\gl_r}$ which is an open subset of
the Hilbert scheme ${\Quot}_{\O_X^\chi}^{r,rM}$, and that
$\M_{\o_r}$ is the good quotient of the smooth
$R^0_{\gl_r}$-scheme $R_{\o_r}$ which represents the functor of
global sections of the quotient by $\o_r$ of the universal bundle
parametrized by $R^0_{\gl_r}$.

\begin{remark}
(i) Note that properness of $\M_{\ssl_r}$ follows from the
construction, since this moduli space is obtained as the good
quotient of the set of semi-stable points of a projective variety.
For arbitrary structure group, this construction does not ensure
the properness of the moduli space (while Ramanathan's original
one did), and we have to use instead semi-stable reduction
theorems for principal $G$-bundles.

\noindent(ii) Existence of moduli spaces for principal bundles has
been since then proved for higher dimensional base varieties and
in arbitrary characteristic (see \cite{gomezetal} and \cite{schmitt}).
\end{remark}

\section{Orthogonal and symplectic bundles}

Let us now specialize the preceding discussion to the
classical groups $\o_r$, $\so_r$ (with $r \geqslant 3$) and
$\Sp_{2r}$. In these cases $\o_r$- and $\Sp_{2r}$-bundles are just
orthogonal and symplectic bundles, and $\so_r$-bundles are
oriented orthogonal bundles:

\begin{definition}
An \textit{orthogonal bundle} is a vector bundle $E$ endowed with
a non-degenerate quadratic form $q \colon E \lra \O_X$ (or,
equivalently, with a symmetric isomorphism $i \colon E \lra
E^\ast$). An \textit{oriented orthogonal bundle} is an orthogonal
bundle $(E,q)$ with an orientation, which comes as a section
$\omega \in H^0(X,\det E)$ of the determinant line bundle of $E$
satisfying $\tilde q (\omega)=1$, where $\tilde q$ is the
quadratic form on $\det E$ deduced from $q$.

A \textit{symplectic} bundle is a vector bundle $E$ endowed with a
non-degenerate symplectic form $\varphi \colon \Lambda^2 E \lra
\O_X$ (or with an antisymmetric isomorphism $E \lra E^\ast$).
\end{definition}

From now on, we concentrate on orthogonal bundles, and generally omit the corresponding statements for symplectic bundles.

\subsection{} \label{ssorth}For these bundles, semi-stability condition \ref{Gss} translates in
a very convenient way: an orthogonal bundle $(E,q)$ is semi-stable
if and only if the underlying vector bundle $E$ is semi-stable.
However, an orthogonal bundle is stable if and only if it splits
as the direct orthogonal sum of some mutually non isomorphic
orthogonal bundles which are stable as vector bundles (see
\cite{ramanan}).

It follows from \ref{unitary} (see also \cite[Theorem
3.18]{ramram}) that if $(E,q)$ is a unitary orthogonal bundle then
$E$ is already a polystable vector bundle. It means that $E$
splits as a direct sum of stable vector bundles. Let us recall two
elementary facts about stable vector bundles: they are simple
bundles, and there are no non zero morphism between non isomorphic
stable vector bundles of the same slope. Hence, the non-degenerate
quadratic structure suggests to write $E$ as
\begin{align*} \num \label{E} E=\bigoplus\limits_{i=1}^{n_1} \Bigl(F^{(1)}_{i}\otimes V^{(1)}_{i}\Bigr) \oplus \bigoplus\limits_{j=1}^{n_2} \Bigl(F^{(2)}_{j} \otimes V^{(2)}_{j}\Bigr) \oplus \bigoplus\limits_{k=1}^{n_3} \Bigl((F^{(3)}_{k} \oplus {F^{(3)}_{k}}{}^\ast)\otimes V^{(3)}_{k}\Bigr),\end{align*}
\noindent where $(F^{(1)}_{i})$ (resp. $(F^{(2)}_{j})$, resp. $(F^{(3)}_{k})$) is a family of mutually non isomorphic orthogonal (resp. symplectic, resp. non isomorphic to their dual nor to that of $F^{(3)}_{k'}$, $k' \neq k$) stable bundles, and $(V^{(1)}_i)$ (resp. $(V^{(2)}_j)$, resp. $(V^{(3)}_k)$) are quadratic (resp. symplectic, resp. equipped with a non-degenerate bilinear form) vector spaces, whose dimension counts the multiplicity of the corresponding stable vector bundle in $E$. Note that the subbundles $F^{(3)}_k \oplus F^{(3)}_k {}^\ast$ have been tacitly endowed with the standard hyperbolic quadratic forms.

\begin{remark}
 We gave in \cite[Remark 1.3 (ii)]{comp} another way to obtain the previous description of a unitary orthogonal bundle, which is in a sense more algebraic (since it avoids the use of the result of Narasimhan and Seshadri), and shows how the unitary bundle associated to a given orthogonal bundle can be defined in terms of isotropic filtrations of the underlying vector bundle (see also \cite{bhosle}).
\end{remark}

\subsection{} According to Ramanathan's result, there exists a moduli space $\M_{\so_r}$ for semi-stable oriented orthogonal vector bundles of rank $r$ on $X$. It is a projective scheme, whose points correspond to unitary oriented orthogonal bundles. It has two connected components $\M_{\so_r}^+$ and $\M_{\so_r}^-$, which are distinguished by the second Stiefel-Whitney class $w_2 \in H^2(X,\Z/2\Z)=\{ \pm 1 \}$.

We also have a moduli scheme $\M_{\o_r}$ for semi-stable
orthogonal bundles. It has several connected components, each of
them corresponding to orthogonal bundles $E$ of a given
topological type, which is determined here by the Stiefel-Whitney
classes $w_i(E) \in H^i(X,\Z/2\Z)$, $i=1,2$ (note that the first
Stiefel-Whitney class $w_1(E)$ is nothing but the determinant
$\det E$).

These moduli varieties have been investigated for hyperelliptic
curves by Ramanan in \cite{ramanan} and Bhosle in \cite{bhosle}.
More recently, for $g \geqslant 2$, Beauville studied in
\cite{arnorth} the linear system associated to the determinant
line bundle on $\M_{\so_r}^{\pm}$ (which is, for $n \neq 4$, the
ample generator of the Picard group of $\M_{\so_r}^\pm$).

We will consider too the moduli scheme $\M_{\Sp_{2r}}$ of semi-stable symplectic bundles on $X$, which is an irreducible projective variety.

\subsection{} Let us introduce now the forgetful morphism

$$\xymatrix@R=10pt@C=50pt{
\M_{\o_r} \ar[r] & \M_{\gl_r}^0 \\
(E,q) \ar@{|->}[r] & E
}$$
\noindent which forgets the quadratic structure, as well as the other forgetful morphisms $\M_{\so_r} \lra \M_{\ssl_r}$ and $\M_{\Sp_{2r}} \lra \M_{\ssl_{2r}}$. It follows from the construction of the moduli schemes that these morphisms are finite (see e.g. \cite{balsesh}). We give in this section a set-theoretic study of these morphisms.

\begin{propo}
 The forgetful morphisms $\M_{\o_r} \lra \M_{\gl_r}^0$ and $\M_{\Sp_{2r}} \lra \M_{\ssl_{2r}}$ are injective.
\end{propo}

\begin{proof}
It is enough to check injectivity on closed points. In view of
\ref{ssorth} we have to prove that any two quadratic structures on
a given polystable vector bundle $E$ define isomorphic orthogonal
bundles, or, in other words, that they differ by a linear
automorphism of $E$. The decomposition (\ref{E}) of $E$ shows that
its automorphism group is isomorphic to
\begin{align*}  \mathrm{Aut}_{\gl_r}(E)= \prod_{i=1}^{n_1} \gl(V^{(1)}_{i}) \times \prod_{j=1}^{n_2} \gl(V^{(2)}_{j}) \times \prod_{k=1}^{n_3} \Bigl(\gl(V^{(3)}_{k}) \times \gl(V^{(3)}_{k}) \Bigr).\end{align*}
Since a quadratic structure on $E$ is nothing but the data of non-degenerate quadratic (resp. symplectic, resp. bilinear) forms on the vector spaces $V_i^{(1)}$ (resp. $V_j^{(2)}$, resp. $V_k^{(3)}$), the conclusion simply follows from the basic fact that any two non-degenerate quadratic (or symplectic) forms on a vector space over an algebraically closed field are equivalent.
\end{proof}

\begin{remark} This proposition is in fact a very particular case of a more general result proved by Grothendieck in \cite{grothendieck} (see also \cite{balaji}). Indeed, this result holds for any vector bundle on any projective variety defined over an algebraically closed field (of characteristic different from $2$).
\end{remark}

\subsection{} Oriented orthogonal bundles behaves differently. For orthogonal bundles of odd rank, $-1$ gives an orthogonal automorphism which exchanges orientation, and it follows that two $\so_r$-structures on a vector bundle of odd rank are automatically equivalent. In even rank, this is no longer true. In fact it already fails for rank $2$ bundles. Indeed, $\so_2$ is just the multiplicative group $\G_m$, which means that $L \oplus L^{-1}$ and $L^{-1} \oplus L$ (endowed with their oriented hyperbolic form) are not isomorphic as $\so_2$-bundles. We can give a precise criterion for unitary orthogonal bundles to admit two non equivalent orientations.

\begin{propo}
 A unitary orthogonal bundle $(E,q)$ has two antecedents via the forgetful morphism
$$\M_{\so_r} \lra \M_{\o_r}$$
\noindent if and only if every orthogonal bundle appearing in the decomposition (\ref{E}) of its underlying vector bundle has even rank.
\end{propo}

\noindent In particular, orthogonal bundles of even rank whose underlying vector bundle is stable have two distinct reductions of structure group to $\so_r$. Such bundles always exist for curves of genus $g \geqslant 2$, and therefore the generic fiber of $\M^\pm_{\so_r} \lra \M_{\o_r}$ has two elements when $r$ is even.

\begin{proof}
 The two orientations defined on a given orthogonal bundle $(E,q)$ give isomorphic $\so_r$-bundles if and only if there is an orthogonal automorphism of $(E,q)$ which exchanges the orientation. It follows from decomposition (\ref{E}) that the isometry group of $E$ is the subgroup of $\Aut_{\gl_r}(E)$ equal to
\begin{align*}\num  \label{Auto} {\Aut}_{\o_r}(E)=\prod\limits_{i=1}^{n_1} {\o}({V^{(1)}_{i}}) \times \prod\limits_{j=1}^{n_2} \Sp(V^{(2)}_{j}) \times \prod\limits_{k=1}^{n_3} \gl(V^{(3)}_{k}),\end{align*}
\noindent where $\gl({V^{(3)}_{k}})$ is identified with its image
in $\gl({V^{(3)}_{k}}) \times \gl({V^{(3)}_{k}})$ by the morphism
$g \mapsto (g, {}^t g^{-1})$. So $E$ admits orthogonal
automorphisms with non trivial determinant if and only if at least
one of the bundles $F^{(1)}_i$ has odd rank.

 Note that this argument is encoded by the exact sequence of non-abelian cohomology associated to the exact sequence $1 \to \so_r \to \o_r \to \mu_2 \to 1$ (see \cite[Remark 1.6]{comp}).
\end{proof}

\subsection{} Before closing this section, we would like to describe precisely what happens in the case of elliptic curves. Moduli spaces of $G$-bundles on an elliptic curve have been described in \cite{laszloell}: if we denote by $\Gamma(T)$ the group of one parameter subgroups of a maximal torus $T \subset G$, the connected component of topologically trivial $G$-bundles is the quotient of $X \otimes_\Z \Gamma(T)$ by the operation of the Weyl group $W_T$. We give here a direct elementary proof of this fact for orthogonal and symplectic bundles.

\begin{propo}
Let $X$ be an elliptic curve, and $l \geqslant 1$. The moduli space $\mathcal{M}_{\so_{2l+1}}^+ $ is isomorphic to $\P ^l$, $\mathcal{M}_{\so_{2l+1}}^- $ to
$\P^{l-1}$,  $\mathcal{M}_{\so_{2l}}^- $ to $\P^{l-2}$
and $ \mathcal{M}_{\so_{2l}}^+ $ to the quotient of $X^l$ by $\left(\Z / 2 \Z \right)^{l-1} \rtimes \mathfrak{S}_l$.
\end{propo}

\begin{proof}
We know by Atiyah's classification that every semi-stable vector bundle of degree zero on the elliptic curve $X$ is $S$-equivalent to a direct sum of line bundles. In particular, if $\kappa_1, \kappa_2$ and $\kappa_3$ are the three non zero line bundles of order $2$, an orthogonal bundle $E$ on $X$ with trivial determinant splits as follows:
\begin{itemize}
\item $\O_X  \oplus \bigoplus\limits_{i=1}^l \left(L_i \oplus L_i^{-1}\right)$ if $\mathrm{rk}(E)=2l+1$ and $w_2(E)=1$,
\item $\kappa_1 \oplus \kappa_2 \oplus \kappa_3 \oplus \bigoplus\limits_{i=1}^{l-1}\left(L_i \oplus L_i^{-1}\right)$ if $\mathrm{rk}(E)=2l+1$ and $w_2(E)=-1$,
\item $\bigoplus\limits_{i=1}^l \left(L_i \oplus L_i^{-1}\right)$ if $\mathrm{rk}(E)=2l$ and $w_2(E)=1$,
\item $\O_X \oplus \kappa_1 \oplus \kappa_2 \oplus \kappa_3 \oplus \bigoplus\limits_{i=1}^{l-2}\left(L_i \oplus L_i^{-1}\right)$ if $\mathrm{rk}(E)=2l$ and $w_2(E)=-1$,
\end{itemize}
\noindent where the $L_i$ are degree $0$ line bundles on $X$. In
all cases but the third one, there is at least one line bundle of
order $2$ which allows us to adjust the determinant of an
orthogonal isomorphism: in these cases, we see that closed points
of the moduli spaces are characterized by collections
$\{M_1,\ldots,M_k\}$ where $M_i \in \{L_i,L_i^{-1}\}$. This gives
the expected isomorphisms, since $X^k/\bigl((\Z/2\Z)^k \rtimes
\mathfrak S_k\bigr)$ is the $k$-th symmetric product of $\P^1$,
which is isomorphic to $\P^k$.

In the remaining case, a generic orthogonal bundle admits two unequivalent orientations, and $\M_{\so_{2l}}^+$ is a quotient of $X^l$ by the action of $(\Z/2\Z)^{l-1} \rtimes \mathfrak S_l$ where $(\Z/2\Z)^{l-1}$ acts on $X \times \cdots \times X$ by transformations $(a_1,\ldots,a_l) \longmapsto (\pm a_1,\ldots,\pm a_l)$ with an even number of minus signs. This finishes the proof of the proposition.

(Of course, a complete proof would consist in defining morphisms from the products of copies of $X$ to the corresponding moduli spaces, and checking that these morphisms induce the above isomorphisms.)
\end{proof}

\begin{remark}\label{elliptic}
(i) In particular, the forgetful morphism $\M_{\so_r}^- \lra
\M_{\ssl_r}$ is always a closed immersion for an elliptic curve,
contrary to what happens in higher genus.

\noindent (ii) When $r$ is even, the moduli space $\M_{\o_r}^\circ$ of topologically trivial orthogonal bundles with trivial determinant is isomorphic to $\P^{\frac{r}{2}}$, and the forgetful morphism $\M_{\so_r}^+ \lra \M_{\o_r}^\circ$ is a $2$-sheeted covering.

\noindent (iii) Of course, the same argument applies to moduli of symplectic bundles and gives isomorphisms between $\M_{\mathbf{Sp}_{2r}}$ and $\P^r$.
\end{remark}

\section{Differential study of the forgetful morphism and quiver representations}
\label{imm}

We have seen that the forgetful morphism $\M_{\o_r} \lra \M_{\gl_r}^0$ is an injective finite morphism. It is natural to ask whether we can say more about this morphism. The answer is given by the main result of \cite{comp}.

\begin{theorem} \label{main}
 The forgetful morphism $\M_{\o_r} \lra \M_{\gl_r}$ is a closed immersion.
\end{theorem}

Of course, the symplectic version of this statement also holds.

\begin{theorem}
 The forgetful morphism $\M_{\Sp_{2r}} \lra \M_{\ssl_{2r}}$ is a closed immersion.
\end{theorem}

\subsection{}  Before going into details, let us give a few remarks about the proof. Since the forgetful morphism is injective and proper, it remains to show that it is everywhere locally a closed immersion, or, equivalently, that it is unramified (see \cite[17.2.6]{EGAIV4} and \cite[8.11.5]{EGAIV3}). Here again, it is enough to consider closed points.

To do this, we use Luna's \'etale slice theorem to get a good
enough understanding of the local structure of the moduli spaces
$\M_{\o_r}$ and $\M_{\gl_r}^0$: we thus obtain \'etale affine
neighbourhoods which appear as good quotients of affine spaces by
the action of some reductive groups. At this point, we have to
understand the corresponding coordinate rings, which are exactly
the invariant rings associated to these actions, and to check that
the ring morphism induced by the forgetful morphism is surjective.
In particular, it is enough to find generating sets for these
invariant rings.

\subsection{} We begin by exhibiting \'etale neighbourhoods for moduli spaces of vector bundles.

\begin{lemma}[{\cite[Theorem 1]{localstr}}] \label{locgl}
 At a polystable vector bundle $E$, the moduli scheme $\M_{\gl_r}^0$ is \'etale locally isomorphic to a neighbourhood of the origin in the good quotient
$$\Ext^1(E,E) \git \Aut_{\gl_r}(E),$$
\noindent where $\Aut_{\gl_r}(E)$ acts on $\Ext^1(E,E)$ by functoriality.
\end{lemma}

\begin{proof}
The construction sketched in \ref{construction} presents
$\M_{\gl_r}^0$ as the good quotient of a smooth open subscheme
$R_{\gl_r}^0$ of the Hilbert scheme $\Quot_{\O_X^\chi}^{r,rM}$ by
the natural action of the reductive group $\Gamma=\gl_\chi$. Let
$q \in R^0_{\gl_r}$ be a point over $E \in \M_{\gl_r}^0$ whose
orbit $\Gamma \cdot q$ is closed, and denote by $N_q$ the normal
space at $q$ to this orbit. We know by Luna's \'etale slice
theorem that there exists a locally closed subscheme $V \subset
R^0_{\gl_r}$ containing $q$ and invariant for the action of the
isotropy group $\Gamma_q \subset \Gamma$ of $q$, together with a
$\Gamma_q$-equivariant morphism $V \lra N_q$ sending $q$ onto $0$,
such that the morphisms
$$V \git \Gamma_q \lra \M_{\gl_r}^0 \ \text{ and } \  V \git \Gamma_q \lra N_q \git \Gamma_q$$
\noindent are \'etale. Deformation theory shows that $N_q$ is
isomorphic to the space of extensions $\Ext^1(E,E)$ of $E$ by
itself, while an easy argument proves that the isotropy group
$\Gamma_q$ is isomorphic to $\Aut_{\gl_r}(E)$.
\end{proof}

\subsection{}
Let us now carry out the same analysis for orthogonal bundles. Recall that, if $P=(E,q)$ is an orthogonal bundle, we denote by $\Ad(P)$ its adjoint bundle $\Ad(P)=P\times^{\o_r} \mathfrak{so}_r$. The symmetric isomorphism $\sigma \colon E \lra E^\ast$ given by the quadratic structure identifies $\Ad(P)$ with the subbundle of $\mathcal End(E)$ consisting of germs of endomorphisms $f$ satisfying $\sigma f + f^\ast \sigma=0$ (or, equivalently, with the vector bundle $\Lambda^2 E^\ast$). The first cohomological space $H^1(X,\Ad(P))$ is thus isomorphic to the space $\Ext^1_\mathrm{asym}(E,E) \subset \Ext^1(E,E)$ of antisymmetric extensions of $E$ by itself.

\begin{lemma}\label{loco}
 At a unitary orthogonal bundle $P=(E,q)$, the moduli scheme is \'etale locally isomorphic to a neighbourhood of the origin in
$$H^1(X,\Ad(P))\git \Aut_{\o_r}(P).$$
\noindent Moreover, the forgetful morphism coincides, through the
different local isomorphisms, to the natural morphism
$$H^1(X,{\Ad}(P)) \git {\Aut}_{\o_r}(P) \lra \Ext^1(E,E) \git \Aut_{\gl_r}(E)$$
induced by the inclusion $H^1(X,\Ad(P)) \subset \Ext^1(E,E)$.
\end{lemma}

\begin{proof}
It follows from \ref{construction} that $\M_{\o_r}$ is the
quotient of a smooth $R_{\gl_r}^0$-scheme $R_{\o_r}$ by the group
$\Gamma$. Hence Luna's theorem applies as well as in the case of
vector bundles: if $q'$ is a point of $R_{\o_r}$ with closed orbit
lying over $P$, $N_{q'}$ the normal space at $q'$ to this orbit,
and $\Gamma_{q'}$ the isotropy group of $q'$, we can find a slice
$V'$ through $q'$ together with a $\Gamma_{q'}$-equivariant
morphism $V'\lra N_{q'}$ giving \'etale morphisms
$$V' \git \Gamma_{q'} \lra \M_{\o_r} \ \text{ and } \  V' \git \Gamma_{q'} \lra N_{q'} \git \Gamma_{q'}.$$
Deformation theory implies that the normal space $N_{q'}$ is
isomorphic to $H^1(X,\Ad(P))=\Ext^1_\mathrm{asym}(E,E)$, and we
can check that the isotropy group $\Gamma_{q'}$ is isomorphic to
$\Aut_{\o_r}(E)$ (we abusively write $\Aut_{\o_r}(E)$ instead of
$\Aut_{\o_r}(P)$).

The second part follows from the fact that the forgetful morphism is the quotient by $\Gamma$ of the structural morphism $R_{\o_r} \lra R_{\gl_r}^0$. We may then choose compatible slices $V$ and $V'$ in order to obtain the following commutative diagram
$$\xymatrix@C=50pt@R=15pt{
& V'\git \Gamma_{q'} \ar[dl] \ar[dr] \ar[dd] & \\
\Ext^1_{\text{asym}}(E,E) \git \mathrm{Aut}_{\o_r}(E) \ar[dd] &  & \M_{\o_r}  \ar[dd]\\
 & V \git \Gamma_q\ar[rd] \ar[ld] & \\
\Ext^1(E,E) \git \mathrm{Aut}_{\gl_r}(E)  &  & \M^0_{\gl_r},
}$$
\noindent which gives the expected identification.
\end{proof}

\subsection{}\label{restr}
We have thus translated the infinitesimal study of the forgetful morphism to a question regarding the morphism
$$\Ext^1_{\text{asym}}(E,E) \git \mathrm{Aut}_{\o_r}(E) \lra \Ext^1(E,E) \git \mathrm{Aut}_{\gl_r}(E).$$
\noindent Theorem \ref{main} is proved if we
show that, for every unitary orthogonal bundle $(E,q)$, this
morphism is unramified at the origin. Now, if we denote by $k[X]$ the coordinate ring of an affine scheme $X$, this morphism corresponds to the restriction morphism
$$k[\Ext^1(E,E)]^{\Aut_{\gl_r}(E)} \lra k[\Ext^1_{\text{asym}}(E,E)]^{\Aut_{\o_r}(E)}$$
between invariant algebras, and it is enough to check that it is a surjective morphism.

\begin{remark}
 On the open locus of $\M_{\o_r}$ consisting of orthogonal bundles with stable underlying vector bundle, Theorem \ref{main} is automatic. Indeed the isotropy groups act trivially, and there is nothing left to prove.
\end{remark}

\subsection{} \label{quiver}We now make Lemmas \ref{locgl} and \ref{loco} more explicit. The polystable vector bundle $E$ can be written as
$$\displaystyle{E = \bigoplus_{i=1}^n F_i \otimes V_i}$$
\noindent where $F_1,\ldots,F_n$ are mutually non isomorphic stable vector bundles, and $V_1,\ldots,V_n$ vector spaces. The space of extensions $\Ext^1(E,E)$ decomposes as
$$\Ext^1(E,E)=\bigoplus_{i,j} \Ext^1(F_i,F_j)\otimes \Hom(V_i,V_j)$$
and the isotropy group is isomorphic to $\Aut_{\gl_r}(E)=\prod_i
\gl(V_i)$. Denote by $d_{ij}$ the dimension of $\Ext^1(F_i,F_j)$,
which is equal to $\rk(F_i) \rk (F_j) (g-1)$ for $i \neq j$, and
to $\rk(F_i)^2 (g-1)+1$ for $i=j$. Thus, if we pick bases for the
extension spaces $\Ext^1(F_i,F_j)$, we may view an extension
$\omega \in \Ext^1(E,E)$ as a collection
$(f_{ij}^k)^{}_{1\leqslant i,j\leqslant n,\ k=1,\ldots,d_{ij}}$ of
morphisms between the vector spaces $V_1,\ldots,V_n$. An element
$g=(g_1,\ldots,g_n)\in \Aut_{\gl_r}(E)$ acts on $(f_{ij}^k)$ by
conjugation:
$$g \cdot (f_{ij}^k) = (g^{}_j f_{ij}^k g^{-1}_i).$$

We recognize here the setting of quiver representations (see \cite{michel}). Indeed, let us consider the quiver $\mathcal Q_E$ whose set of vertices is defined by
$$(\mathcal Q_E)_0=\{s_1,\ldots,s_n\},$$
\noindent these vertices being connected by $d_{ij}$ arrows from $s_i$ to $s_j$, and define the dimension vector $\alpha \in \N^{(\mathcal Q_E)_0}$ by $\alpha_i=\dim V_i$. The preceding discussion shows that $\Ext^1(E,E)$ is exactly the representation space $R(\mathcal Q_E,\alpha)$ of the quiver $\mathcal Q_E$ for the dimension vector $\alpha$, and that the action of $\Aut_{\gl_r}(E)$ on $\Ext^1(E,E)$ is nothing but the usual action of the group $\gl(\alpha)=\prod_i \gl_{\alpha_i}$ on $R(\mathcal Q_E,\alpha)$:

\begin{lemma}\label{quivergl}
 The $\mathrm{Aut}_{\gl_r}(E)$--module $\mathrm{Ext}^1(E,E)$ is isomorphic to the $\prod \gl_{\alpha_i}$--module $R(\mathcal Q_E,\alpha)$. In particular, it only depends (up to $\mathrm{Aut}_{\gl_r}(E)$-isomorphism) on the ranks and multiplicities of the stable subbundles $F_1,\ldots,F_n$ of $E$.
\end{lemma}

\subsection{} \label{extasym} Suppose now that $(E,q)$ is a unitary orthogonal bundle. Following (\ref{E}), we can write $E$ as the direct sum
\begin{align*}  E=\bigoplus\limits_{i=1}^{n_1} \Bigl(F^{(1)}_{i}\otimes V^{(1)}_{i}\Bigr) \oplus \bigoplus\limits_{j=1}^{n_2} \Bigl(F^{(2)}_{j} \otimes V^{(2)}_{j}\Bigr) \oplus \bigoplus\limits_{k=1}^{n_3} \Bigl((F^{(3)}_{k} \oplus {F^{(3)}_{k}}{}^\ast)\otimes V^{(3)}_{k}\Bigr).\end{align*}
\noindent Let us put $E^{(1)}_i= F^{(1)}_i \otimes V^{(1)}_i$, $E^{(2)}_j= F^{(2)}_j \otimes V^{(2)}_j$ and $E^{(3)}_k= (F^{(3)}_k \oplus F^{(3)}_k{}^\ast)\otimes V^{(3)}_k$. They all have an orthogonal structure $ \sigma^{(a)}_l \colon E^{(a)}_l \buildrel\sim\over\to E^{(a)}_l{}^\ast$ induced by that of $E$.

The space $\mathrm{Ext}^1(E,E)$ splits into the direct sum of all extension spaces $\mathrm{Ext}^1(E^{(k)}_i,E^{(l)}_j)$. An element $\omega=\sum \omega^{(k,l)}_{i,j} \in \mathrm{Ext}^1(E,E) \simeq \bigoplus \mathrm{Ext}^1(E^{(k)}_i,E^{(l)}_j)$ belongs to $\Ext^1_\mathrm{asym}(E,E)$ if and only if $\omega^{(k,k)}_{i,i} \in \Ext^1_\mathrm{asym}(E^{(k)}_i,E^{(k)}_i) \subset \mathrm{Ext}^1(E^{(k)}_i,E^{(k)}_i)$ for all $i$ and $k$, and $\sigma^{(l)}_j \omega^{(k,l)}_{i,j} + {\omega^{(l,k)}_{j,i}}{}^\ast \sigma^{(k)}_i=0$ for all $(i,k) \neq (j,l)$. So, identifying $\mathrm{Ext}^1(E_i^{(k)},E_j^{(l)})$ with its image in $\mathrm{Ext}^1(E_i^{(k)},E_j^{(l)}) \oplus \mathrm{Ext}^1(E_j^{(l)},E_i^{(k)})$ by the application $\omega_{i,j}^{(k,l)} \mapsto \omega_{i,j}^{(k.l)} - {\sigma_i^{(k)}}{}^{-1}{\omega_{i,j}^{(k,l)}}{}^\ast \sigma_j^{(l)}$, it appears that $\Ext^1_\mathrm{asym}(E,E)$ is equal to the subspace of $\mathrm{Ext}^1(E,E)$ defined as
\begin{align*}
 \bigoplus\limits_{k} \left(\bigoplus\limits_i  \Ext^1_\mathrm{asym}(E^{(k)}_i,E^{(k)}_i) \oplus \bigoplus\limits_{i<j} \mathrm{Ext}^1(E^{(k)}_i,E^{(k)}_j) \right) \oplus
\bigoplus\limits_{k<l} \bigoplus\limits_{i,j}  \mathrm{Ext}^1(E^{(k)}_i,E^{(l)}_j).
\end{align*}
\noindent Moreover, we can check that the diagonal summands involved in this decomposition are respectively isomorphic to:
\begin{align*}   \Ext^1_{\mathrm{asym}}(E^{(1)}_i,E^{(1)}_i)= \Bigl( H^1(X,\mathrm S^2 F^{(1)}_i{}^\ast) \otimes \sok({V^{(1)}_{i}})\Bigr)  \oplus \Bigl(H^1(X,\Lambda^2 F^{(1)}_i{}^\ast) \otimes \mathrm{S}^2 {V^{(1)}_{i}}{}^\ast\Bigr),
\end{align*}
\begin{align*}
   \Ext^1_{\mathrm{asym}}(E^{(2)}_j,E^{(2)}_j)=\Bigl( H^1(X,\Lambda^2 F^{(2)}_j{}^\ast) \otimes \spk({V^{(2)}_{j}})\Bigr) \oplus \Bigl( H^1(X,\mathrm S^2F^{(2)}_j{}^\ast) \otimes {\Lambda}\!^2 {V^{(2)}_{j}}{}^\ast\Bigr),
\end{align*}
\begin{multline*}
  \Ext^1_{\mathrm{asym}}(E^{(3)}_k,E^{(3)}_k)= \Bigl(\mathrm{Ext}^1(F^{(3)}_{k},F^{(3)}_{k}) \otimes \glk({V^{(3)}_{k}})\Bigr) \oplus \\ \Bigl(H^1(X,\mathrm S^2 F^{(3)}_k{}^\ast) \otimes {\Lambda}^2V^{(3)}_{k}{}^\ast\Bigr) \oplus  \Bigl(H^1(X, \Lambda^2 F^{(3)}_k{}^\ast) \otimes \mathrm{S}^2{V^{(3)}_{k}}{}^\ast\Bigr) \oplus \\ \Bigl(H^1(X, \mathrm S^2 F^{(3)}_k) \otimes {\Lambda}^2V^{(3)}_{k}\Bigr) \oplus \Bigl(H^1(X,\Lambda^2 F^{(3)}_k) \otimes \mathrm{S}^2{V^{(3)}_{k}}\Bigr),
\end{multline*}
\noindent where $\mathrm{Ext}^1(F^{(3)}_{k},F^{(3)}_{k})$ has been identified with its image in $\mathrm{Ext}^1(F^{(3)}_{k},F^{(3)}_{k}) \oplus \mathrm{Ext}^1({F^{(3)}_{k}}{}^\ast,{F^{(3)}_{k}}{}^\ast)$ by the map $\omega \mapsto \omega-\omega^\ast$.

The isometry group
\begin{align*} {\Aut}_{\o_r}(E)=\prod\limits_{i=1}^{n_1} {\o}({V^{(1)}_{i}}) \times \prod\limits_{j=1}^{n_2} \Sp(V^{(2)}_{j}) \times \prod\limits_{k=1}^{n_3} \gl(V^{(3)}_{k}),\end{align*}
\noindent (see (\ref{Auto})) naturally acts on $\Ext_\mathrm{asym}^1(E,E)$ by conjugation.

This laborious description of the  $\Aut_{\o_r}(E)$-module $\Ext^1_\mathrm{asym}(E,E)$ has the following consequence:

\begin{lemma}\label{ind}
 The morphism
$$\Ext^1_\mathrm{asym}(E,E) \git \Aut_{\o_r}(E) \lra \Ext^1(E,E) \git \Aut_{\gl_r}(E)$$
\noindent only depends, up to isomorphisms, on the ranks and multiplicities of the stable bundles $F^{(a)}_l$ appearing in the decomposition (\ref{E}) associated to the orthogonal bundle $E$.
\end{lemma}

\subsection{Case of the trivial bundle}\label{trivial}
In order to clarify a bit this description before proving the main
result of this section (as well as to give an idea of this proof),
it seems useful to consider the case of the trivial orthogonal
bundle $E=\O_X \otimes k^r$. The space of extensions $\Ext^1(E,E)$
is then identified with the space $\Mat_r(k)^g$ of $g$-tuples of
$r \times r$ matrices, and $\Ext^1_\mathrm{asym}(E,E)$ with the
subspace $\Mat_{r}^\mathrm{asym}(k)^g$ of $g$-tuples of
antisymmetric matrices. The isotropy groups
$\Aut_{\gl_r}(E)=\gl_r$ and $\Aut_{\o_r}(E)=\o_r$ act diagonally
by conjugation.

As we have seen in \ref{restr}, the forgetful morphism is unramified at the trivial bundle if the restriction morphism
$$k[\Mat_r(k)^g]^{\gl_r} \lra k[\Mat_r^\mathrm{asym}(k)^g]^{\o_r}$$
\noindent is surjective. These invariant algebras have been
described in \cite{PrAdv}. The algebra $k[\Mat_r(k)^g]^{\gl_r}$ is
generated by traces of products $(M_1,\ldots,M_g) \longmapsto
\tr(M_{i_1} \cdots M_{i_l})$ (with $l\leqslant r^2$), while
$k[\Mat_r(k)^{g}]^{\o_r}$ is generated by functions
$(M_{1},\ldots,M_g) \longmapsto \tr(A_{i_1} \cdots A_{i_l})$,
where $A_{i_k} \in \{M_{i_k},{}^t M_{i_k}\}$. The restriction of
such a function to the subspace $\Mat_r^\mathrm{asym}(k)^g$ is
clearly the restriction of a $\gl_r$-invariant function on
$\Mat_r(k)^g$. Since the restriction map $k[\Mat_r(k)^g]^{\o_r}
\to k[\Mat_r^\mathrm{asym}(k)^g]^{\o_r}$ is surjective, it proves
that the forgetful morphism is indeed unramified at the trivial
bundle.

\begin{proof}[Proof of Theorem \ref{main}]
 We have to prove that the forgetful morphism is unramified. The decomposition (\ref{E}) of a unitary orthogonal bundle allows us to define the \textit{slice-type stratification} of $\M_{\o_r}$. The locally closed strata consist of all unitary orthogonal bundles $E$ having a given isometry group $\Aut_{\o_r}(E)$. Lemma \ref{ind}, together with Lemma \ref{loco}, shows that the sheaf of relative differential
 $$\Omega^1_{\M_{\o_r}/\M_{\gl_r}^0}$$
\noindent has constant rank on each stratum. It is thus enough to show that this sheaf vanishes on the closed ones.

Since orthogonal summands $F^{(3)}_k \oplus F^{(3)}_k{}^\ast$ (with $F^{(3)}_k \not\simeq F^{(3)}_k{}^\ast$) or  $F^{(2)}_j \otimes V^{(2)}_j$ (where $F^{(2)}_j$ is a symplectic bundle) specialize to the trivial orthogonal bundle, a closed stratum must consist of unitary orthogonal bundles which split as
$$E = \bigoplus_{i=1}^{n_1} F^{(1)}_i \otimes V^{(1)}_i$$
\noindent where $F^{(1)}_1,\ldots,F^{(1)}_{n_1}$ are mutually non isomorphic orthogonal bundles whose underlying vector bundles are stable, and $V^{(1)}_1,\ldots,V^{(1)}_{n_1}$ some quadratic spaces.

Let now $E=\bigoplus_i F_i \otimes V_i$ be such an orthogonal bundle. We claim that the restriction morphism
$$k[\Ext^1(E,E)]^{\Aut_{\gl_r}(E)} \lra k[\Ext^1_{\text{asym}}(E,E)]^{\Aut_{\o_r}(E)}$$
\noindent is surjective. This means that the forgetful morphism is unramified on the closed strata, which finishes the proof of the Theorem.

Let $\mathcal Q_E$ be the quiver defined in \ref{quiver}.
According to Lemma \ref{quivergl}, $k[\Ext^1(E,E)]^{\Aut_{\gl_r}(E)} $ is isomorphic to the coordinate ring $k[R(\mathcal Q_E,\alpha)]^{\gl(\alpha)}$ of the quotient variety $R(\mathcal Q_E,\alpha)\git \gl(\alpha)$ which parametrizes isomorphism classes of semisimple representations of $\mathcal Q_E$ with dimension $\alpha$. This invariant algebra has been described in \cite{LBP} (in the characteristic $0$ case). In particular, it is generated by traces along oriented cycles in $\mathcal Q_E$ (of length $\leqslant (\sum \alpha_i)^2$), that is by functions
$$(f_a)^{}_a \longmapsto \mathrm{tr}(f_{a_l} \cdots f_{a_1})$$
where $a_1 \cdots a_l$ is an oriented cycle in the quiver $\mathcal Q_E$.

On the other hand, since the inclusion $\Ext^1_\mathrm{asym}(E,E) \lra \Ext^1(E,E)$ is equivariant for the action of the isometry group
$$\Aut_{\o_r}(E)=\prod_i \o(V_i),$$
\noindent the restriction morphism $k[\Ext^1(E,E)]^{\Aut_{\o_r}(E)} \lra k[\Ext^1_{\text{asym}}(E,E)]^{\Aut_{\o_r}(E)}$ is surjective. The next proposition provides us with a set of generators for
$$k[\Ext^1(E,E)]^{\Aut_{\o_r}(E)} \simeq k[R(\mathcal Q_E,\alpha)]^{\prod \o_{\alpha_i}}.$$
\noindent If $\widetilde{\mathcal Q_E}$ is the quiver deduced from $\mathcal Q_E$ by adding one new arrow $a^\ast \colon v' \to v$ for any arrow $a \colon v \to v'$, it tells us that $k[R({\mathcal Q_E},\alpha)]^{\prod \o_{\alpha_i}}$ is generated by functions
$$ (f_a)^{}_a \longmapsto \mathrm{tr}(f_{\tilde{a}_l} \cdots f_{\tilde{a}_1})$$
where $\tilde{a}_1 \cdots \tilde{a}_l$ is an oriented cycle in the quiver $\widetilde{\mathcal Q_E}$, and $f_{\tilde{a}_i}$ is equal to $f_{a_i}$ or its adjoint according to whether $\tilde{a}_i$ is $a_i$ or ${a_i}^\ast$.

The space $\Ext^1_\mathrm{asym}(E,E) \subset \Ext^1(E,E)$ has been described in \ref{extasym}. It identifies with a subspace of $R(\mathcal Q_E,\alpha)$ made of representations having the following property: if $f_a \colon V_v \lra V_{v'}$ is the linear morphism associated to an arrow $a \colon v \to v'$, then its adjoint morphism $f_a^\ast \colon V_{v'}^\ast \lra V_v^\ast$ is, up to the sign, the linear morphism associated to an arrow $a \colon v' \to v$. It implies that the restrictions to $\Ext^1_\mathrm{asym}(E,E)$ of the preceding functions are also restrictions of  $\Aut_{\gl_r}(E)$-invariant functions on $\Ext^1(E,E)$, whence our claim.
\end{proof}

Let us now state and prove the result about $\prod \o_{\alpha_i}$-invariant functions used in the proof. Let $Q$ by a quiver with $n$ vertices, and $\alpha \in \N^n$ a dimension vector. Consider the group $\o(\alpha)=\prod \o_{\alpha_i}$. As a subgroup of $\gl(\alpha)$, it acts by conjugation on the representation space $R(Q,\alpha)$. Let $\widetilde Q$ be the quiver deduced from $Q$ by adding one new arrow $a^\ast \colon v' \to v$ for any arrow $a \colon v \to v'$.

\begin{propo}\label{inv}
 The algebra $k[R(Q,\alpha)]^{\o(\alpha)}$ of polynomial invariants for the action of $\o(\alpha)$ on the representation space $R(Q,\alpha)$ is generated by traces along oriented cycles in the associated quiver $\widetilde Q$. These are functions
$$(f_a)^{}_a \mapsto \mathrm{tr}(f_{\tilde{a}_l} \cdots f_{\tilde{a}_1})$$
where $\tilde{a}_1 \cdots \tilde{a}_l$ is an oriented cycle in the quiver $\widetilde Q$, and $f_{\tilde{a}_i}$ is equal to $f_{a_i}$ or its adjoint according to whether $\tilde{a}_i$ is $a_i$ or ${a_i}^\ast$.
\end{propo}

It may be rephrased as follows. First note that any representation
of $Q$ can be extended to a representation of $\widetilde Q$ by
associating to a new arrow $a^\ast$ the adjoint of the linear map
corresponding to $a$. This defines a natural map $R(Q,\alpha) \lra
R(\widetilde Q,\alpha)$, and the
proposition just means that the restriction morphism
$k[R(\widetilde Q,\alpha)]^{\gl(\alpha)} \lra
k[R(Q,\alpha)]^{\o(\alpha)}$ is onto.

This result is a special case of \cite[Theorem 2.3.3]{comp}, and
follows (exactly as (\textit{loc. cit.})) from an adaptation of
the proof given in \cite{LBP} to describe the invariant ring
$k[R(Q,\alpha)]^{\gl(\alpha)}$. This special case is technically
much easier. Indeed, we had to consider in \cite{comp} algebras
with antimorphisms of order $4$, while antiinvolutions are enough
here. We present here a quite detailed proof, but warmly refer the
reader to the original exposition \cite{LBP}.

We first need a lemma about algebras with trace and antiinvolution. Recall that a $k$-algebra with trace is a $k$-algebra $R$ together with a linear map $\tr \colon R \to R$ satisfying the identities $ \tr(a) b=b\, \tr(a)$, $\tr(ab)=\tr(ba)$ and $\tr(\tr(a) b)=\tr(a) \tr(b)$ for all $a, b \in R$. A {$k$-algebra with trace and antiinvolution} is an algebra with trace $R$ endowed with an antiinvolution $\iota \colon R \to R$. The algebra $\Mat_N(B)$ of $N \times N$ matrices with coefficients in a commutative ring $B$ will be equipped with its usual trace together with the adjunction map  $\tau \colon M \in \Mat_N(B) \mapsto {}^t M$.

If $R$ is a $k$-algebra with trace and antiinvolution $\iota$, we can consider the functor $\widetilde X_{R,N}$ (from commutative $k$-algebras to sets) of $N$-dimensional trace preserving representations of $R$ commuting with the antiinvolutions:
$$\widetilde X_{R,N}(B) = \{ f \in \mathrm{Hom}_k(R,\Mat_N(B))\,|\ f\circ\tr=\tr\circ f, f \circ \iota=\tau \circ f \}.$$
\noindent We claim that this functor is representable by a commutative ring $\widetilde A$. Indeed, we know from \cite[2.2]{deconcinietal} that the functor $X_{R,N}$ of trace preserving representations of $R$ is representable by a ring $A$. If $j \colon R \lra \Mat_N(A)$ is the corresponding universal morphism, there is a unique involution $t$ of $A$ such that $\Mat_N(t) j=\tau j \iota$, and the quotient $\widetilde A$ of $A$ by this involution represents $\widetilde X_{R,N}$. We still denote by $\widetilde X_{R,N}$ the affine scheme $\Spec \widetilde A$.

In particular, we have a universal morphism $\tilde j \colon R \lra \Mat_N(\widetilde A)$. The conjugation action of $\o_N$ on $\Mat_N(\widetilde A)$ induces a right action on $\widetilde A$: indeed, every $g$ in $\o_N$ defines an automorphism $\bar{g}$ of $\widetilde A$ such that $\Mat(\bar{g})\tilde j=C(g)\tilde j$, where we denote by $C(g)$ the conjugation by $g$. We consider the action of $\o_N$ on $\Mat_N(\widetilde A)$ defined by $g \cdot M=C(g) \Mat(\bar{g})^{-1}(M)$ for $g \in \o_N$ and $M \in \Mat_N(\widetilde A)$. The universal morphism $\tilde j$ then maps $R$ into the algebra $\Mat_N(\widetilde A)^{\o_N}$ of $\o_N$-equivariant morphisms from $\X_{R,N}={\Spec\widetilde A}$ to $\Mat_N(k)$ (see \cite{PrJofalg} or \cite[1.2]{deconcinietal}).

The main result of \cite{PrJofalg} can be easily adapted to this situation (see also \cite[\S 12]{berele}):

\begin{lemma}\label{pr}
 Let $R$ be a $k$-algebra with trace and antiinvolution. Then the universal morphism $\tilde j$ is a surjective morphism $R \lra \Mat_N(\widetilde A)^{\o_N}$.
\end{lemma}

\begin{proof}  Following \cite{PrJofalg} we begin by proving this when $R$ is a free algebra with trace and antiinvolution built on the generators $\{x_s\}_{s \in \Sigma}$. In this case one can check that $\Mat_N(\widetilde A)^{\o_N}$ is the algebra of all $\o_N$-equivariant polynomial maps from $\Mat_N(k)^\Sigma$ to $\Mat_N(k)$, and our assertion immediately follows from the description of this algebra given in \cite[7.2]{PrAdv} (which comes as a direct consequence of the result recalled in \ref{trivial} about generators for the invariant algebra $k[\Mat_N(k)^\Sigma]^{\o_N}$).

In the general case, we write $R$ as the quotient of a free algebra with trace and antiinvolution $T$ by an ideal $I$. If $\widetilde A _T$ is the universal ring associated to $T$, we know that the two-sided ideal in $\Mat_N(\widetilde A _T)$ generated by the image of $I$ must be equal to $\Mat_N(J)$ for some ideal $J$ in $\widetilde A_T$. The universal ring for $R$ is then the quotient $\widetilde A_T / J$. The conclusion follows from the linear reductivity of $\o_r$, which ensures that $\Mat_N(\widetilde A_T)^{\o_r} \lra \Mat_N(\widetilde A_T/J)^{\o_r}$ is onto. Note that this last argument makes essential use of the characteristic $0$ assumption.
\end{proof}

Let us go back to the quiver $Q$. The associated quiver
$\widetilde Q$ carries a natural involution $\sigma$ that fixes
vertices and exchanges arrows $a$ and $a^\ast$. Let $R$ (resp.
$\widetilde R$) be the algebra obtained from the path algebra of
the opposite quiver $Q^{\mathrm{op}}$ (resp. $\widetilde
Q^{\mathrm{op}}$) by adding traces. The involution $\sigma \colon
\widetilde Q \to \widetilde Q$ induces an antiinvolution $\iota$
of $\widetilde R$ such that representations of $\widetilde R$
commuting with $\tau$ and $\iota$ correspond bijectively to
representations of $R$. In other words, $\sigma$ gives an
involution of the space $R(\widetilde Q,\alpha)$ such that
$R(Q,\alpha)$ is isomorphic to the subspace of $R(\Q,\alpha)$
consisting of all representations which preserve the preceding
involutions. The proof of \ref{inv} relies on a precise
description of this space as a subspace of $\X_{\widetilde
R,N}(k)$, where $N=\sum \alpha_i$.

\begin{proof}[Proof of Proposition \ref{inv}] We follow closely \cite[\S 3]{LBP}. Consider the subalgebra $\widetilde S_n \subset \widetilde R$ generated by the orthogonal idempotents $e_1,\ldots,e_n$ corresponding to the different vertices of $\widetilde Q$. The antiinvolution $\iota$ is trivial on this subalgebra, and the scheme $\widetilde X_{\widetilde S_n,N}$ is the disjoint union $\widetilde X_{\widetilde S_n,N}=\bigcup_\delta \widetilde X_\delta$, where $\delta \in \N^n$ ranges over the set of all dimension vectors such that $\sum \delta_i=N$, of the homogeneous varieties
$$\widetilde X_{\delta}=\o_N/\prod_i \o_{\delta_i}.$$
\noindent  This induces a decomposition $\widetilde X_{\widetilde R,N}=\bigcup \varpi^{-1}\widetilde X_\delta$ (where $\varpi \colon \widetilde X_{\widetilde R,N} \lra \widetilde X_{\widetilde S_n,N}$ is the morphism induced by the inclusion $\widetilde S_n \subset \widetilde R$), or, equivalently, a decomposition $\widetilde A=\prod \widetilde A_\delta$ of the coordinate ring $\widetilde A$ of $\widetilde X_{\widetilde R,N}$.

We focus on the component $\widetilde X_{\widetilde
R,\alpha}=\varpi^{-1}\widetilde X_\alpha$ corresponding to the
dimension vector $\alpha$. Let us write the identity matrix
$\mathrm{id_N}$ as the sum $\sum u_i$ of orthogonal idempotents
$u_1,\ldots,u_n$ associated to the decomposition $k^N= \bigoplus
k^{\alpha_i}$, and let $p$ be the point in $\widetilde X_\alpha$
defined by the representation of $\widetilde S_n$ sending $e_i$ to
$u_i$. The fiber $\varpi^{-1}(p)$ (which represents the subfunctor
of $\widetilde X_{\widetilde R,N}$ consisting of representations
sending $e_i$ to $u_i$) naturally carries an action of the
centralizer in $\o_N$ of the idempotents $u_i$, which is
isomorphic to $\o(\alpha)=\prod \o_{\alpha_i}$. Moreover, this
fiber can be identified with the subspace of $R(\widetilde
Q,\alpha)$ consisting of involutions preserving representations of
$\widetilde Q$, which is itself isomorphic to $R(Q,\alpha)$.

Since $\varpi$ is $\o_N$-equivariant and $\widetilde X_\alpha$ is homogeneous, the invariant ring $\Mat_N(k[R(Q,\alpha)])^{\o(\alpha)}$ of $\o_\alpha$-equivariant maps from $\varpi^{-1}(p)$ to $\Mat_N(k)$ is exactly the ring $\Mat_N(\widetilde A_\alpha)^{\o_N}$ of $\o_N$-equivariant maps from $\widetilde X_{\widetilde R,\alpha}$ to $\Mat_N(k)$. But, since $\o_N$ acts separately on each factor of $\Mat_N(\widetilde A)=\prod \Mat_N(\widetilde A_\delta)$, it follows from Lemma \ref{pr} that $\tilde j$ gives a surjective morphism from $\widetilde R$ onto $\Mat_N(\widetilde A_\alpha)^{\o_N} \simeq \Mat_N(k[R(Q,\alpha)])^{\o(\alpha)}$. The expected description of $k[R(Q,\alpha)]^{\o(\alpha)}$ follows by taking traces.
\end{proof}

\begin{remark}
(i) We have treated in \cite{comp} the more general following problem: let $Q$ stand for a quiver with $n=n_1+n_2+n_3+2n_4$ vertices
$$r_1,\ldots,r_{n_1},s_1, \ldots, s_{n_2}, t_1,\ldots,t_{n_3}, u_1,u_1^\ast,\ldots,u_{n_4},u_{n_4}^\ast,$$
\noindent and $\alpha \in \N^n \simeq \N^{n_1} \times \N^{n_2} \times \N^{n_3} \times (\N \times \N)^{n_4}$ be an admissible dimension vector (by which we mean a vector such that  $\alpha_{t_k}$ is even and $\alpha_{u_l}=\alpha_{u_l^\ast}$). We define $\Gamma(\alpha)$ to be the group
$$\Gamma(\alpha)=\prod_{i=1}^{n_1} \gl_{\alpha_{r_i}} \times \prod_{j=1}^{n_2} \o_{\alpha_{s_j}} \times \prod_{k=1}^{n_3} \Sp_{\alpha_{t_k}} \times \prod_{l=1}^{n_4} \gl_{\alpha_{u_l}},$$
\noindent which is actually thought of here as a subgroup of $\gl(\alpha)$ via the inclusions $P \in \gl_{\alpha_{u_l}} \mapsto (P,{}^tP^{-1}) \in \gl_{\alpha_{u_l}} \times \gl_{\alpha_{u_l^\ast}}$ for $l=1,\ldots,n_4$. We give a generating set for the invariant algebra $k[R(Q,\alpha)]^{\Gamma(\alpha)}$.

Let $\widetilde Q$ be the quiver deduced from $Q$ as follows. We add $n_1$ new vertices $r_1^\ast, \ldots,r_{n_1}^\ast$, and consider the involution $\sigma$ on the set of vertices $\widetilde Q_0$ which fixes the $s_j$ and $t_k$ and exchanges $u_l$ with $u_l^\ast$, and $r_i$ with $r_i^\ast$. We now add a new arrow $a^\ast \colon \sigma (v') \to \sigma(v)$ for any arrow $a \colon v \to v'$. Note that in this new quiver two vertices $r_i$ and $r_{i'}^\ast$ are never connected by a single arrow.

\begin{theoremnonnum}
The invariant algebra $k[R(Q,\alpha)]^{\Gamma(\alpha)}$ is generated by traces along cycles in the associated quiver $\widetilde Q$, that is by functions
$$(f_a)^{}_a \mapsto \mathrm{tr}(f_{\tilde{a}_l} \cdots f_{\tilde{a}_1})$$
where $\tilde{a}_1 \cdots \tilde{a}_l$ is an oriented cycle in the quiver $\widetilde Q$, and $f_{\tilde{a}_i}$ is equal to $f_{a_i}$ or its adjoint according to whether $\tilde{a}_i$ is $a_i$ or ${a_i}^\ast$.
\end{theoremnonnum}

It directly implies (together with \ref{extasym}) that $\Ext^1_\mathrm{asym}(E,E) \git \Aut_{\o_r}(E) \lra \Ext^1(E,E)\git \Aut_{\gl_r}(E)$ is a closed immersion for every unitary orthogonal bundle $E$.

\vspace{0.2cm}

\noindent (ii) The proof given in \cite{LBP} and its adaptation in Proposition \ref{inv} are only valid in characteristic $0$. However, these results remain true in arbitrary characteristic. Characteristic free proofs can be found in \cite{donkin} for the $\gl(\alpha)$-action on $R(Q,\alpha)$ and in \cite{lopatin} for the $\Gamma(\alpha)$-action. They rely on the notions of good filtrations and good pairs.
\end{remark}

\subsection{} We close this section with a few comments about the forgetful morphism $\M_{\so_r} \lra \M_{\ssl_r}$ for genus $\geqslant 2$ curves (see Remark \ref{elliptic} for elliptic curves). Of course its study reduces to that of $\M_{\so_r} \lra \M_{\o_r}$. Lemma \ref{loco} and its variant for $\so_r$-bundles easily show that, when $r$ is odd, this morphism is an isomorphism onto its image. When $r$ is even, it is a $2$-sheeted covering of its image. Indeed, an orthogonal bundle $(E,q)$ has two antecedents if and only if $\Aut_{\so_r}(E)=\Aut_{\o_r}(E)$. Then, at a point $(E,q,\omega)$ defined by any of its two reductions to $\so_r$, the forgetful morphism is a local isomorphism (in the \'etale topology).

\begin{propo}
  Let $X$ be a curve of genus $\geqslant 2$. Then the forgetful morphism $\M_{\so_r} \lra \M_{\ssl_r}$ is a closed immersion when $r$ is odd, and a $2$-sheeted covering of its image when $r$ is even.
\end{propo}

\begin{remark}
 When $r$ is even, the infinitesimal behaviour of $\M_{\so_r} \lra \M_{\o_r}$ remains difficult to describe explicitly. In the case of the trivial bundle $E=\O \otimes V$ with $V$ a quadratic vector space of even dimension, we have to understand the action of ${\Aut}_{\so_r}(E) \simeq \so_r$ on $\Ext_\mathrm{asym}^1(E,E) \simeq H^1(X,\O_X) \otimes \sok(V)$. The computation has been carried out in \cite{aslaksen}, and provides a set of generators for the $k[\Ext_\mathrm{asym}^1(E,E)]^{\mathrm{Aut}_{\o_r}(P)}$-algebra $k[\Ext_\mathrm{asym}^1(E,E)]^{\mathrm{Aut}_{\so_r}(P)}$ in terms of \textit{polarized pfaffians}.
\end{remark}

\section{Singular locus of $\M_{\so_r}$ and $\M_{\Sp_{2r}}$}\label{sing}

In this section, $X$ is a curve of genus $g \geqslant 2$.

\subsection{} Narasimhan and Ramanan have described in \cite{NR} the singular locus of the moduli space $\M_{\gl_r}^d$ of vector bundles of rank $r$ and degree $d$: they have shown that it is exactly the closed subset of strictly semi-stable vector bundles, except when $X$ is a curve of genus $2$, $r=2$ and $d$ even (in which case $\M_{\gl_2}^d$ is a smooth variety). Note that the fact that stable vector bundles $E$ define non singular points in $\M_{\gl_r}^d$ is trivial. Indeed, we know that $\M_{\gl_r}^d$ is \'etale locally isomorphic to $\Ext^1(E,E) \git \Aut_{\gl_r}(E)$, which is smooth when $E$ is stable for the obvious reason that the isotropy group $\Aut_{\gl_r}(E)=\G_m$ then acts trivially on $\Ext^1(E,E)$.

\subsection{} For arbitrary reductive algebraic groups $G$, the relevant notion is that of \textit{regularly stable bundle}: a regularly stable $G$-bundle is a stable $G$-bundle $E$ such that $\mathrm{Aut}_G(E)=Z(G)$. The same argument shows that the smooth locus of $\M_G$ contains the open subset of regularly stable bundles.

We check here that, when $G=\so_r$ and $\o_r$, this inclusion is
in fact an equality, except in two special cases (which are not
surprising in view of the particular case occurring in \cite{NR}).
Note that regularly stable oriented orthogonal bundles are stable
orthogonal bundles whose underlying vector bundle is either stable
or, when $r$ is even, the direct sum of two different stable
bundles of odd rank, while regularly stable orthogonal bundles are
just orthogonal bundles with stable underlying vector bundle.

\begin{theorem}
 The smooth locus of $\M_{\so_r}$ (resp. $\M_{\o_r}$) is precisely the open set consisting of regularly stable $\so_r$-bundles (resp. $\o_r$-bundles), except when $g=2$ and $r=3$ or $4$.
\end{theorem}

\begin{proof}
The proof relies on the precise description of the closed points
of $\M_{\so_r}$. Let $U$ be the set of points in $\M_{\so_r}$
which correspond to those oriented orthogonal bundles
$(E,q,\omega)$ which are either regularly stable, or an orthogonal
sum $E=E_1 \oplus E_2$ of two different stable vector bundles, or
a symplectic sum $E=F \oplus F$ of two copies of a regularly
stable symplectic bundle, or an hyperbolic sum $E=F \oplus F^\ast$
where $F$ is a stable vector bundle with $F \not \simeq F^\ast$.
It is an open dense subset of $\M_{\so_r}$. Moreover, if we denote
by $\M^\rs_{\so_r}$ the locus of regularly stable bundles, we see
that $U \setminus \M_{\so_r}^\rs$ is open and dense in $\M_{\so_r}
\setminus \M^\rs_{\so_r}$ for $r \geqslant 4$ (when $r=3$ we need
to enlarge $U$ by adding the bundles $\O_X \oplus L \oplus L^{-1}$
for degree $0$ line bundles $L$ with $L^2 \not \simeq \O_X$).

It is thus enough to check that the singular locus of $U$ is
exactly $U \setminus \M_{\so_r}^\rs$. The proof of Lemma
\ref{loco} shows that $\M_{\so_r}$ is \'etale locally isomorphic
at a point defined by a unitary bundle $P$ to the good quotient
$H^1(X,\Ad(P))\git \Aut_{\so_r}(P)$.

At a point defined by an oriented orthogonal bundle $E=E_1 \oplus
E_2$ with $E_1$ and $E_2$ two different stable vector bundles (of
even rank if $r$ is even), $\M_{\so_r}$ is locally isomorphic to
an \'etale neighbourhood of the origin in the quotient of
$$H^1(X,\Lambda^2 E_1) \oplus H^1(X,\Lambda^2 E_2) \oplus \Ext^1(E_1,E_2)$$
\noindent by the action of $\mathrm{Aut}_{\so_r}(E) / Z(\so_r)
\simeq \mu_2$ (where $-1 \in \mu_2$ acts by $(1,1,-1)$).
Chevalley's theorem implies that this quotient cannot be smooth,
since $\Ext^1(E_1,E_2)$ must have dimension at least $2$.

At a point $E=F \oplus F^\ast$ with $F$ a stable vector bundle non isomorphic to its dual, $\M_{\so_r}$ is locally isomorphic to the quotient of
$$H^1(X,\Lambda^2 F^\ast) \oplus H^1(X, \Lambda^2 F) \oplus \Ext^1(F,F)$$
\noindent by the action of $\mathrm{Aut}_{\so_r}(E) \simeq
\mathrm{Aut}(F) \simeq \mathbf G_m$, where $\lambda \in \mathbf
G_m$ acts by $(\lambda^{-2},\lambda^2,1)$. We easily see that its
multiplicity at the origin is equal to $\genfrac{(}{)}{0pt}{}{2
(d-1)}{d-1}$ where $d$ is the dimension of $H^1(X,\Lambda^2
F^\ast)$ (see \cite[3.3.4 (ii)]{comp}). This multiplicity cannot
be equal to $1$.

Finally, an \'etale neighbourhood of a point $E=F \otimes V$
defined by a regularly stable symplectic bundle $F$ and a
symplectic vector space $V$ of dimension $2$ is given by an
\'etale neighbourhood of the origin in the quotient of
$$\left(H^1(X,\mathrm{S}^2 F) \otimes \Lambda^2 V \right) \oplus \left( H^1(X,\Lambda^2 F) \otimes \mathrm{S}^2 V\right)$$
\noindent by the action of $\mathrm{Aut}_{\so_r}(E) = \Sp(V)
\simeq \Sp_2$. It follows from the classification of all coregular
representations of almost simple connected complex algebraic
groups given in \cite{schwarz} that this quotient cannot be smooth
unless $\dim H^0(C,\Lambda^2 F) \leqslant 2$, which cannot happen
but for a rank $2$ symplectic bundle $F$ on a curve of genus $2$.

This concludes the proof of the Theorem for $r \geqslant 4$. In
rank $3$ we have also to consider the points $E=\O_X \oplus L
\oplus L^{-1}$ where $L$ is a line bundle of degree $0$ whose square
$L^2$ is not trivial. The automorphism $\lambda \in
\mathrm{Aut}_{\so_3}(E)\simeq \mathrm{Aut} (L)=\mathbf G_m$ then acts
on
$$\Ext^1(L,L) \oplus \Ext^1(\O_X,L) \oplus \Ext^1(\O_X,L^{-1}) $$
\noindent by $(1,\lambda,\lambda^{-1})$. We easily see that the
quotient is smooth if and only if $\Ext^1(\O_X,L)$ has dimension
$1$, which happens exactly when $X$ has genus $2$.
\end{proof}

\begin{remark}
If $g=2$ and $r=3$, the same techniques can be used to describe
the singular locus of $\M_{\so_3}$. Since the sum of two copies of
the adjoint representation of $\so_3$ is coregular, the trivial
bundle $\O_X \oplus \O_X \oplus \O_X$ defines a smooth point in
$\M_{\so_3}$ (even if this particular point is often called ``the
worst point''). So the singular locus is exactly the closure of
the set of orthogonal bundles of the form $\eta \oplus F$ where
$\eta \not\simeq \O_X$ is a line bundle of order $2$ and $F$ a
rank $2$ orthogonal bundle with $\det(F)=\eta$ which is stable as
a vector bundle.

If $r=4$, the smooth locus is exactly the union of
$\M_{\so_4}^\rs$ and the locally closed subset corresponding to
orthogonal bundles which are symplectic sums $F^{(2)} \otimes
V^{(2)}$ of two copies of a stable (symplectic) bundle.
\end{remark}

We can of course prove in the same way the following result for moduli
spaces $\M_{\Sp_{2r}}$ of semi-stable symplectic bundles on $X$ of
rank $2r \geqslant 4$:

\begin{theorem}
The smooth locus of $\M_{\Sp_{2r}}$ is precisely the open set consisting of regularly stable bundles.
\end{theorem}

\backmatter

\nocite{}
\bibliographystyle{smfplain}
\bibliography{SCgrenoble}

\providecommand{\bysame}{\leavevmode ---\ }
\providecommand{\og}{``}
\providecommand{\fg}{''}
\providecommand{\smfandname}{\&}
\providecommand{\smfedsname}{\'eds.}
\providecommand{\smfedname}{\'ed.}
\providecommand{\smfmastersthesisname}{M\'emoire}
\providecommand{\smfphdthesisname}{Th\`ese}
\begin{thebibliography}{10}

\bibitem{aslaksen}
{\scshape H.~Aslaksen, E.-C. Tan {\normalfont \smfandname} C.-B. Zhu} -- {\og
  Invariant theory of special orthogonal groups\fg}, \emph{Pacific J. Math.}
  \textbf{168} (1995), no.~2, p.~207--215.

\bibitem{balaji}
{\scshape V.~Balaji} -- {\og Lectures on principal bundles\fg}, in \emph{Moduli
  Spaces and Vector Bundles}, London Math. Soc. Lecture Note Ser., vol. 359,
  Cambridge Univ. Press, Cambridge, 2009.

\bibitem{balsesh}
{\scshape V.~Balaji {\normalfont \smfandname} C.~S. Seshadri} -- {\og
  Semistable principal bundles--{I} (characteristic zero)\fg}, \emph{J.
  Algebra} \textbf{258} (2002), no.~1, p.~321--347.

\bibitem{beauville2}
{\scshape A.~Beauville} -- {\og Fibr\'es de rang {$2$} sur une courbe, fibr\'e
  d\'eterminant et fonctions th\^eta\fg}, \emph{Bull. Soc. Math. France}
  \textbf{116} (1988), no.~4, p.~431--448.

\bibitem{arnorth}
\bysame , {\og Orthogonal bundles on curves and theta functions\fg}, \emph{Ann.
  Inst. Fourier (Grenoble)} \textbf{56} (2006), no.~5, p.~1405--1418.

\bibitem{BLS}
{\scshape A.~Beauville, Y.~Laszlo {\normalfont \smfandname} C.~Sorger} -- {\og
  The {P}icard group of the moduli of {$G$}-bundles on a curve\fg},
  \emph{Compositio Math.} \textbf{112} (1998), no.~2, p.~183--216.

\bibitem{berele}
{\scshape A.~Berele} -- {\og Matrices with involution and invariant theory\fg},
  \emph{J. Algebra} \textbf{135} (1990), no.~1, p.~139--164.

\bibitem{bhosle}
{\scshape U.~N. Bhosle} -- {\og Moduli of orthogonal and spin bundles over
  hyperelliptic curves\fg}, \emph{Compositio Math.} \textbf{51} (1984), no.~1,
  p.~15--40.

\bibitem{michel}
{\scshape M.~Brion} -- {\og Representations of quivers\fg}, in \emph{Geometric
  methods in representation theory, Summer school, Grenoble, 2008}.

\bibitem{deconcinietal}
{\scshape C.~De~Concini, C.~Procesi, N.~Reshetikhin {\normalfont \smfandname}
  M.~Rosso} -- {\og Hopf algebras with trace and representations\fg},
  \emph{Invent. Math.} \textbf{161} (2005), no.~1, p.~1--44.

\bibitem{donkin}
{\scshape S.~Donkin} -- {\og Polynomial of invariants of representations of
  quivers\fg}, \emph{Comment. Math. Helvetici} \textbf{69} (1994), p.~137--141.

\bibitem{gomezetal}
{\scshape T.~G\'omez, A.~Langer, A.~H.~W. Schmitt {\normalfont \smfandname}
  I.~Sols} -- {\og Moduli spaces for principal bundles in arbitrary
  characteristic\fg}, \emph{Adv. Math.} \textbf{219} (2008), no.~4,
  p.~1177--1245.

\bibitem{grothendieck}
{\scshape A.~Grothendieck} -- {\og Sur la classification des fibr\'es
  holomorphes sur la sph\`ere de {R}iemann\fg}, \emph{Amer. J. Math.}
  \textbf{79} (1957), p.~121--138.

\bibitem{EGAIV3}
\bysame , {\og \'{E}l\'ements de g\'eom\'etrie alg\'ebrique, {IV}. \'{E}tude
  locale des sch\'emas et des morphismes de sch\'emas {III}\fg}, \emph{Inst.
  Hautes \'Etudes Sci. Publ. Math.} (1966), no.~28.

\bibitem{EGAIV4}
\bysame , {\og \'{E}l\'ements de g\'eom\'etrie alg\'ebrique, {IV}. \'{E}tude
  locale des sch\'emas et des morphismes de sch\'emas {IV}\fg}, \emph{Inst.
  Hautes \'Etudes Sci. Publ. Math.} (1967), no.~32.

\bibitem{localstr}
{\scshape Y.~Laszlo} -- {\og Local structure of the moduli space of vector
  bundles over curves\fg}, \emph{Comment. Math. Helv.} \textbf{71} (1996),
  no.~3, p.~373--401.

\bibitem{laszloell}
\bysame , {\og About {$G$}-bundles over elliptic curves\fg}, \emph{Ann. Inst.
  Fourier (Grenoble)} \textbf{48} (1998), no.~2, p.~413--424.

\bibitem{LBP}
{\scshape L.~Le~Bruyn {\normalfont \smfandname} C.~Procesi} -- {\og Semisimple
  representations of quivers\fg}, \emph{Trans. Amer. Math. Soc.} \textbf{317}
  (1990), no.~2, p.~585--598.

\bibitem{lopatin}
{\scshape A.~A. Lopatin} -- {\og Invariants of quivers under the action of
  classical groups\fg}, \emph{J. Algebra} \textbf{321} (2009), no.~4,
  p.~1079--1106.

\bibitem{mumfordICM}
{\scshape D.~Mumford} -- {\og Projective invariants of projective structures
  and applications\fg}, in \emph{Proc. {I}nternat. {C}ongr. {M}athematicians
  ({S}tockholm, 1962)}, Inst. Mittag-Leffler, Djursholm, 1963, p.~526--530.

\bibitem{NR}
{\scshape M.~S. Narasimhan {\normalfont \smfandname} S.~Ramanan} -- {\og Moduli
  of vector bundles on a compact {R}iemann surface\fg}, \emph{Ann. of Math.
  (2)} \textbf{89} (1969), p.~14--51.

\bibitem{NarasimhanSeshadri}
{\scshape M.~S. Narasimhan {\normalfont \smfandname} C.~S. Seshadri} -- {\og
  Stable and unitary vector bundles on a compact {R}iemann surface\fg},
  \emph{Ann. of Math. (2)} \textbf{82} (1965), p.~540--567.

\bibitem{PrAdv}
{\scshape C.~Procesi} -- {\og The invariant theory of {$n\times n$}
  matrices\fg}, \emph{Adv. Math.} \textbf{19} (1976), no.~3, p.~306--381.

\bibitem{PrJofalg}
\bysame , {\og A formal inverse to the {C}ayley-{H}amilton theorem\fg},
  \emph{J. Algebra} \textbf{107} (1987), no.~1, p.~63--74.

\bibitem{ramanan}
{\scshape S.~Ramanan} -- {\og Orthogonal and spin bundles over hyperelliptic
  curves\fg}, in \emph{Geometry and analysis}, Indian Acad. Sci., Bangalore,
  1980, p.~151--166.

\bibitem{ramram}
{\scshape S.~Ramanan {\normalfont \smfandname} A.~Ramanathan} -- {\og Some
  remarks on the instability flag\fg}, \emph{T\^ohoku Math. J. (2)} \textbf{36}
  (1984), no.~2, p.~269--291.

\bibitem{ramanathan75}
{\scshape A.~Ramanathan} -- {\og Stable principal bundles on a compact
  {R}iemann surface\fg}, \emph{Math. Ann.} \textbf{213} (1975), no.~5,
  p.~129--152.

\bibitem{ramanathanthesis}
\bysame , {\og Moduli for principal bundles over algebraic curves. {I} and
  {II}\fg}, \emph{Proc. Indian Acad. Sci. Math. Sci.} \textbf{106} (1996),
  p.~301--328, 421--449.

\bibitem{schmitt}
{\scshape A.~H.~W. Schmitt} -- {\og Moduli spaces for principal bundles\fg}, in
  \emph{Moduli spaces and vector bundles}, London Math. Soc. Lecture Note Ser.,
  vol. 359, Cambridge Univ. Press, Cambridge, 2009, p.~388--423.

\bibitem{schwarz}
{\scshape G.~W. Schwarz} -- {\og Representations of simple {L}ie groups with
  regular rings of invariants\fg}, \emph{Invent. Math.} \textbf{49} (1978),
  no.~2, p.~167--191.

\bibitem{comp}
{\scshape O.~Serman} -- {\og Moduli spaces of orthogonal and symplectic bundles
  over an algebraic curve\fg}, \emph{Compos. Math.} \textbf{144} (2008), no.~3,
  p.~721--733.

\bibitem{serre}
{\scshape J.-P. Serre} -- {\og Espaces fibr\'es alg\'ebriques\fg}, S\'em. C.
  Chevalley {\textbf{3}} (1958/59), Expos\'e no.1, 37 p.

\bibitem{LPV}
{\scshape J.-L. Verdier {\normalfont \smfandname} J.~Le~Potier} (\smfedsname)
  -- \emph{Module des fibr\'es stables sur les courbes alg\'ebriques}, Progress
  in Mathematics, vol.~54, Birkh\"auser Boston Inc., Boston, MA, 1985, Papers
  from the conference held at the {\'E}cole Normale Sup{\'e}rieure, Paris,
  1983.

\bibitem{weil}
{\scshape A.~Weil} -- {\og G\'en\'eralisation des fonctions ab\'eliennes\fg},
  \emph{J. de Math. Pures et App.} \textbf{17} (1938), p.~47--87.

\end{thebibliography}

\end{document}